\RequirePackage{fix-cm}
\documentclass[12pt,letterpaper]{article}

\usepackage[margin=1in]{geometry}
\usepackage{amsthm,amsmath,amsfonts,amssymb,mathtools,comment}
\usepackage[authoryear,round]{natbib}
\usepackage{xcolor}
\usepackage[colorlinks=true,linkcolor=blue,citecolor=blue,urlcolor=blue]{hyperref}
\AtBeginDocument{\lefthyphenmin=2\righthyphenmin=3}
\numberwithin{equation}{section}
\theoremstyle{plain}
\newtheorem{theorem}{Theorem}[section]
\newtheorem{proposition}[theorem]{Proposition}
\newtheorem{lemma}[theorem]{Lemma}
\newtheorem{corollary}[theorem]{Corollary}
\newtheorem{remark}{Remark}

\newcommand{\R}{\mathbb R}
\newcommand{\cC}{\mathcal C}
\newcommand{\cL}{\mathcal L}
\newcommand{\cU}{\mathcal U}
\newcommand{\cS}{\mathcal S}
\newcommand{\argmin}{\operatorname*{arg\,min}}
\newcommand{\TV}{\operatorname{TV}}
\newcommand{\lo}{\underline{\widehat\theta}}
\newcommand{\hi}{\overline{\widehat\theta}}

\title{Pointwise Minmax Representations for Total Variation Denoising on Graphs}
\author{Sabyasachi~Chatterjee and Deep~Ghoshal\\[0.5em]
Department of Statistics, University of Illinois Urbana-Champaign.}
\date{}

\begin{document}
\maketitle

\begin{abstract}
Total variation denoising (TVD) is a fundamental technique for image
denoising, signal processing and nonparametric regression, with a long history in statistics and applied mathematics. In this article, we establish exact
pointwise minmax/maxmin representations for its fitted values on
arbitrary finite graphs, extending recent one-dimensional results. Under squared loss, each fitted
value is a minmax of corrected local averages over a pair of nested
connected vertex sets containing the target vertex. The correction is determined by a signed boundary count:
edges leaving the inner set contribute $+1$ if they
remain within the outer set and $-1$ otherwise. This signed boundary count unifies the correction terms
that arise in one dimension and gives the representation
a clean graph-theoretic form. To the best of our knowledge,
this general minmax representation is new even for
two-dimensional TVD, despite its central role in image denoising
and decades of study across multiple scientific communities.

An unexpected connection emerges between these graph TVD
representations and the classical theory of isotonic regression.
We show that Brunk's classical minmax formulas on every finite
directed acyclic graph arise as special cases of our graph TVD
representation, through an exact reduction of isotonic
regression to graph TVD. The reduction passes through
asymmetric TVD and nearly isotonic regression, which agrees
exactly with isotonic regression once the penalty exceeds
a finite threshold.

We further establish pointwise representations under
arbitrary finite separable convex losses, assuming coercivity
of the full objective. Neither differentiability nor strict
convexity of the losses is required. When solutions are nonunique, the
possible fitted values at each coordinate form a compact
interval. The smallest and largest possible fitted values admit maxmin
and minmax representations, respectively, in terms of scalar
local optimization problems and are attained simultaneously
at all coordinates by two extremal solutions. As an application, we obtain new
exact representations for quantile TVD on arbitrary graphs. 

Altogether, our results provide a far-reaching generalization of
Brunk's classical formulas, unifying TVD, nearly isotonic regression,
and isotonic regression within a single pointwise minmax/maxmin
framework for general separable convex losses.
\end{abstract}


\clearpage
 
\section{Introduction and Main Results}\label{sec:introduction}
Total variation denoising (TVD) is a classical nonparametric regression and
denoising method that has been extensively studied in statistics, applied
mathematics, and signal processing. It was introduced in the seminal work of
\cite{rudin1992nonlinear}, which proposed total variation regularization for
two-dimensional image denoising. In one dimension, given observations $y_1,\ldots,y_n$, the
TVD estimator, for $\lambda\geq0$, is
\begin{equation}\label{eq:chain-tvd}
\begin{aligned}
\widehat\theta^{(\lambda)}
&=\argmin_{\theta\in\R^n}\Biggl\{
\frac12\sum_{i=1}^n(y_i-\theta_i)^2\\
&\hspace{25mm}+\lambda\sum_{i=1}^{n-1}|\theta_{i+1}-\theta_i|\Biggr\}.
\end{aligned}
\end{equation}
The total variation (TV) penalty promotes sparsity in the first-order differences
of the fitted signal, and consequently the estimator $\widehat\theta^{(\lambda)}$
is piecewise constant. This estimator is also
known as the fused lasso
\citep{tibshirani2005sparsity} and is closely related to locally adaptive
regression splines \citep{mammen1997locally}.

In the setting of image denoising, the corresponding
TV penalty sums absolute differences across horizontally and
vertically adjacent pixels. If $G=(V,E)$ is the $2D$ grid graph, this penalty is
\begin{equation}\label{eq:graph-tv}
\TV_G(\theta)=\sum_{\{u,v\}\in E}|\theta_u-\theta_v|.
\end{equation}
This is the so-called anisotropic form of
the total variation penalty. The definition in \eqref{eq:graph-tv} extends naturally to any finite
undirected graph, yielding graph TVD, also known as the graph fused lasso.
This estimator encourages piecewise constant signals on connected regions
of the graph; see \citet{tibshirani2011solution} for its formulation within generalized lasso, and
\citet{wang2016trend,hutter2016optimal,padilla2018dfs}
for statistical and computational developments.

The nondifferentiable TV penalty couples neighboring coordinates, making the fitted values of $\widehat\theta^{(\lambda)}$ difficult to characterize explicitly. Despite the long history of TVD, an exact pointwise representation was established only recently in the univariate setting by \citet{chatterjee2026minmax}. This representation is explicit and expresses each fitted value as a min--max, or equivalently a max--min, of simple functions of local averages over nested intervals containing the target point. It also provides a direct route to pointwise statistical analysis, which is difficult to carry out using the global optimization formulation alone. 

In subsequent
work, \citet{ghoshal2026exact} established a corresponding representation
for quantile TVD, replacing the averages by suitable local order
statistics. Quantile loss is neither differentiable nor strictly convex,
and the fitted TVD solution need not be unique. Their formulae characterize
the endpoints of the interval of possible fitted values at each
coordinate.

These two results raise three natural questions, which we answer in this
paper.
\begin{enumerate}
\item \emph{Does a similar representation hold for two-dimensional TVD?}
The results of~\cite{chatterjee2026minmax} and~\cite{ghoshal2026exact} only hold for the one-dimensional case. The main contribution of this article is to show that 
an exact representation holds on every finite
graph, including the $2D$ lattice/grid graph. In one dimension, the representation takes a min--max over nested intervals of local averages with correction terms determined by the position of the inner interval within the outer interval. Our graph representation gives a remarkably clean and natural generalization: intervals become connected vertex sets, and the correction is determined by a count of edges leaving the inner set $I$, with opposite signs according to whether they remain inside or leave the outer set $J$.

\item \emph{How is this representation related to the classical isotonic
minmax formula?} 
The classical min--max formulas for isotonic regression are a cornerstone
of shape-restricted inference: they express a globally constrained estimator through local
averages and provide a direct route to its pointwise statistical
analysis~\cite{robertson1988order,zhang2002risk,chatterjee2015risk}. We show that these formulas
are instances of a more general pointwise representation for graph TVD.
In particular, the classical formula over intersections of upper and
lower sets on a DAG emerges from our formula over nested connected
sets. The bridge is asymmetric TVD, which includes nearly isotonic
regression~\cite{tibshirani2011nearly}, which further agrees exactly with isotonic
regression once the penalty exceeds a finite threshold. This places order-constrained estimation and total
variation regularization within a common min--max framework.

\item \emph{Does the representation extend to a general convex loss?}
The quantile result suggests that such minmax representations can hold for TVD with other (possibly nondifferentiable) losses beyond just the squared loss. In this article, we establish a minmax/maxmin representation
for arbitrary finite, separable convex losses. This general formula recovers both
the squared-loss formula of \citet{chatterjee2026minmax} and the quantile
formula of \citet{ghoshal2026exact}.
\end{enumerate}

Taken together, these results establish a unified pointwise min--max characterization of graph TVD, quantile TVD, nearly-isotonic regression, and isotonic regression. In the next section, we give a precise statement of our pointwise min-max representation for graph TVD.

\subsection{Squared-loss TVD on a general graph}\label{subsec:squared}

Let $G=(V,E)$ be a finite undirected graph with $V\ne\varnothing$. Given
$y\in\R^V$, define the total variation denoising estimator (for the graph $G$) as follows:
\begin{equation}\label{eq:graph-tvd}
\widehat\theta^{(\lambda)}
=\argmin_{\theta\in\R^V}
\left\{\frac12\sum_{u\in V}(y_u-\theta_u)^2
+\lambda\TV_G(\theta)\right\}
\end{equation}
where 
\begin{equation}\label{eq:graph-tv-recalled}
\TV_G(\theta)=\sum_{\{u,v\}\in E}|\theta_u-\theta_v|.
\end{equation}
The objective is strictly convex and coercive, so the minimizer exists
and is unique. Write $\cC$ for the collection of nonempty vertex sets
whose induced subgraphs are connected, and let
$\cC_i=\{J\in\cC:i\in J\}$. For a nonempty set $I$, write
$\bar y_I=|I|^{-1}\sum_{u\in I}y_u$. For disjoint sets $A,B\subseteq V$,
let $|E(A,B)|$ denote the number of edges with one endpoint in $A$ and the other in
$B.$ 

For any pair of nested subsets $\varnothing\ne I\subseteq J\subseteq V$, let us define
\begin{equation}\label{eq:boundary}
C_{I,J}=|E(I,J\setminus I)|-|E(I,J^c)|
\end{equation}
where $J^c=V\setminus J$. Every edge leaving $I$ contributes either
$+1$, if it remains inside $J$, or $-1$, if it leaves $J$ as well.
Thus $C_{I,J}$ is a signed count of the boundary of $I$ (number of edges leaving $I$) with respect to the outer set $J.$ We are now ready to state the main result of this paper.

\begin{theorem}[Pointwise representation for graph TVD]\label{thm:squared}
For every finite undirected graph $G = (V,E)$, every $y\in\R^V$, every
$\lambda\geq0$, and every $i\in V$,

\begin{align}
\widehat\theta_i^{(\lambda)}
&=\max_{J\in\cC_i}
  \min_{\substack{I\in\cC_i\\I\subseteq J}}
  \left\{\bar y_I+\frac{\lambda C_{I,J}}{|I|}\right\}\label{eq:square-maxmin}\\
&=\min_{J\in\cC_i}
  \max_{\substack{I\in\cC_i\\I\subseteq J}}
  \left\{\bar y_I-\frac{\lambda C_{I,J}}{|I|}\right\}.\label{eq:square-minmax}
\end{align}

\end{theorem}

\begin{remark}
In the above result, the connectedness restrictions on \(I\) and \(J\) can be dropped. Both identities remain valid when \(I\) and \(J\) range over arbitrary vertex sets satisfying \(i\in I\subseteq J\).
\end{remark}

\subsubsection{How to read the min--max representation}

Fix a vertex \(i\). Choose an outer connected set \(J\) containing \(i\),
and then an inner connected set \(I\) such that \(i\in I\subseteq J\).
The pair \((I,J)\) gives a corrected local average: the ordinary average
\(\bar y_I\), adjusted by \(\lambda C_{I,J}/|I|\). In the max--min formula
this correction is added; in the min--max formula it is subtracted.
The inner optimization compares these corrected averages for a fixed
\(J\), and the outer optimization then selects among the possible
choices of \(J\).

Let us consider the toy case when $n = 2$, where the graph has two vertices and one edge, and see what the minmax representation gives us. 
Write \(\bar y=(y_1+y_2)/2\). To find the fitted value at vertex \(1\),
there are two choices of outer set. If \(J=\{1\}\), necessarily
\(I=\{1\}\), and \(C_{I,J}=-1\). If \(J=\{1,2\}\), the two choices are
\(I=\{1\}\), for which \(C_{I,J}=1\), and \(I=\{1,2\}\), for which
\(C_{I,J}=0\). Thus the min--max identity reads
\[
\begin{aligned}
\widehat\theta_1^{(\lambda)}
&=\min\bigl\{y_1+\lambda,\max\{y_1-\lambda,\bar y\}\bigr\}\\
&=\operatorname{proj}_{[y_1-\lambda,\,y_1+\lambda]}(\bar y).
\end{aligned}
\]

By symmetry, the corresponding expression for vertex \(2\) projects
\(\bar y\) onto \([y_2-\lambda,y_2+\lambda]\). Writing \(d=y_2-y_1\)
and evaluating these projections gives
\[
\begin{gathered}
(\widehat\theta_1^{(\lambda)},\widehat\theta_2^{(\lambda)})\\
=\begin{cases}
(\bar y,\bar y), & |d|\leq 2\lambda,\\[2mm]
(y_1+\lambda,\;y_2-\lambda), & d>2\lambda,\\[2mm]
(y_1-\lambda,\;y_2+\lambda), & d<-2\lambda.
\end{cases}
\end{gathered}
\]
Hence observations separated by at most \(2\lambda\) are fused at
their average. Otherwise, the two fitted values remain distinct,
with each observation shifted toward the other by \(\lambda\).

\subsubsection{Comparison with the 1D Result}
On a one-dimensional chain graph, connected vertex sets are intervals.
An interval has at most two boundary edges, and \eqref{eq:boundary}
recovers the correction factors in \citet{chatterjee2026minmax}.
For example, when \(J\) avoids the endpoints \(\{1,n\}\),
\(C_{I,J}=2\) if \(I\) lies strictly inside \(J\),
\(C_{I,J}=-2\) if \(I=J\), and \(C_{I,J}=0\) if \(I\)
shares exactly one endpoint with \(J\).

In \citet{chatterjee2026minmax}, these correction factors were
specified case by case, without a unifying formula.
Equation~\eqref{eq:boundary} reveals their underlying structure:
\(C_{I,J}\) is a signed count of the boundary edges of \(I\),
with each edge connecting \(I\) to \(J\setminus I\) contributing
\(1\), and each edge connecting \(I\) to \(J^{c}\)
contributing \(-1\).
This interpretation is purely graph theoretic and makes no reference to the linear ordering of
the $1D$ chain graph. A central contribution of this paper is to identify this
graph-theoretic structure and establish that the same signed
boundary count yields an exact pointwise min--max representation
on every finite graph.

\subsubsection{Implications for two-dimensional TVD}

On a two-dimensional grid, the representation involves arbitrary
connected vertex sets containing the vertex of interest. Unlike
connected sets on the $1D$ chain, these sets need not have a simple
geometric form: they can have irregular boundaries and are not
restricted to rectangles. Thus the local averages appearing in the
formula range over a much richer collection of regions than in
one dimension. 

The representation does not by itself provide a practical method
for computing the estimator. On a two-dimensional grid, the number
of admissible pairs \(I\subseteq J\) can be prohibitively large,
making direct enumeration infeasible. Even in one dimension,
computation is better handled by specialized algorithms, such as
the \(O(n)\) dynamic programming algorithm of
\citet{johnson2013dynamic} or the fast direct algorithm of
\citet{condat2013direct}, without evaluating the min--max formula.

We anticipate that the representation will instead be useful for
studying the pointwise statistical properties of two-dimensional
TVD. It expresses each fitted value through maxima and minima of explicit local averages with boundary corrections, making the roles of the observations
and the geometry of the graph directly visible. In particular,
fixing an outer set \(J\) in the two representations gives upper
and lower bounds on the fitted value. Such bounds offer a starting
point for analyzing pointwise estimation error in terms of local
signal variation, noise averages, and local geometry of the graph.

\subsubsection{A preview of the proof}
\label{subsubsec:proof-intuition}

Why does the representation involve an outer set $J$, an inner set
$I\subseteq J$, and a min--max over these choices? An instructive starting point is the following elementary
identity, valid for every vector $v\in\R^V$:
\begin{equation}\label{eq:vector-minmax}
v_i
=
\min_{\substack{J\subseteq V\\i\in J}}
\max_{\substack{I\subseteq J\\i\in I}}
\bar v_I,
\qquad
\bar v_I=\frac{1}{|I|}\sum_{u\in I}v_u.
\end{equation}
To see this, for any $J$ containing $i$, choose
\[
I=J\cap\{u\in V:v_u\geq v_i\}.
\]
Then $\bar v_I\geq v_i$, so the inner maximum is at least $v_i$
for every $J$. Conversely, choose
\[
J=\{u\in V:v_u\leq v_i\}.
\]
For this particular $J$, every $I\subseteq J$ containing $i$
satisfies $\bar v_I\leq v_i$, proving the reverse inequality.

The exact same choices of level sets also prove the TVD representation. Write
$\widehat\theta=\widehat\theta^{(\lambda)}$.
For every nonempty $I\subseteq V$, the optimality conditions give
\[
\overline{\widehat\theta}_I
=
\bar y_I+\frac{\mathcal F(I,I^c)}{|I|},
\]
where $\mathcal F(I,I^c)$ can be interpreted as a \textit{signed flow} across
the boundary of $I$, going from $I$ to $I^{c}$. Apply the preceding level-set choices with
$v=\widehat\theta$. For the upper bound, the same inner set $I$
inside any $J$ turns out to satisfy
$\mathcal F(I,I^c)\leq-\lambda C_{I,J}$. For the lower bound,
the same outer set $J$ gives the opposite inequality,
$\mathcal F(I,I^c)\geq-\lambda C_{I,J}$, for every
$I\subseteq J$ containing $i$. These flow bounds allow the fitted
averages in the elementary argument to be replaced by the corrected
data averages $\bar y_I-\lambda C_{I,J}/|I|$, with both inequalities
going in the required directions.

We feel that this level-set viewpoint makes the mechanics of the
proof more transparent than in the earlier one-dimensional
arguments of \citet{chatterjee2026minmax} and
\citet{ghoshal2026exact}. The full proof is given in
Section~\ref{sec:proof-squared}.

\subsection{Isotonic Regression is a Special Case of TVD}

It is well known that isotonic regression admits pointwise minmax/maxmin representations. We first recall the classical formula. Let \(D=(V,A)\) be a finite
directed acyclic graph (DAG), and write \(u\preceq v\) if \(u=v\) or
there is a directed path from \(u\) to \(v\). The isotonic cone is
\begin{equation}\label{eq:isotonic-cone-intro}
\mathcal K_D
=\{\theta\in\R^V:\theta_u\leq\theta_v
  \text{ for all }(u,v)\in A\}.
\end{equation}
A lower set contains every predecessor of each of its vertices;
an upper set contains every successor. Equivalently, no directed
edge enters a lower set from outside, and no directed edge leaves
an upper set. Let \(\cL_i\) and \(\cU_i\) denote the lower and upper
sets containing \(i\). The least-squares
isotonic regression estimator is defined by
\[
\widehat\theta^{\mathrm{iso}}
=
\operatorname*{argmin}_{\theta\in\mathbb R^V:\,
\theta_u\leq\theta_v\ \text{for all }(u,v)\in A}
\frac12\sum_{u\in V}(y_u-\theta_u)^2.
\]
Its classical pointwise
representation is
\begin{equation}\label{eq:classical-isotonic-intro}
\max_{U\in\cU_i}\min_{L\in\cL_i}\bar y_{L\cap U}
=
\widehat\theta_i^{\mathrm{iso}}
=
\min_{L\in\cL_i}\max_{U\in\cU_i}\bar y_{L\cap U}.
\end{equation}
This formula plays a major role in the classical theory of order-restricted
inference; see \citet{brunk1955maximum},
\citet{barlow1972statistical}, and
\citet[Theorem~1.4.4]{robertson1988order}.

The resemblance to our TVD representation is immediate, but the
connection is less apparent. In one dimension, the isotonic formula
optimizes over the two endpoints of a single interval, whereas
the TVD formula optimizes over two nested intervals. Moreover,
the isotonic formula involves only local averages, without any
penalty correction. On a general graph, the distinction is between
intersections of upper and lower sets for isotonic regression
and pairs of arbitrary nested connected sets for TVD.
In these respects, the isotonic formula appears simpler.
Nevertheless, as we will argue, it is a special case of our TVD representation.

Write \(\operatorname{ISO}(y,D)=\widehat\theta^{\mathrm{iso}}\) for the isotonic estimator and
\(\operatorname{TVD}(y,G,\lambda)\) for symmetric TVD. The following
proposition makes their connection explicit.

\begin{proposition}[Isotonic regression as TVD]
\label{prop:isotonic-as-tvd}
Let \(D=(V,A)\) be a finite DAG with \(n=|V|\geq1\), and let
\(\widetilde G\) be its underlying undirected graph. Define
\[
d_u=\deg_D^+(u)-\deg_D^-(u),
\qquad
\widetilde y_u^{(\lambda)}
=
y_u-\frac{\lambda}{2}d_u,
\]
where \(\deg_D^+(u)\) and \(\deg_D^-(u)\) denote the out-degree
and in-degree of \(u\). Then, for every
\[
\lambda>\lambda_0
:=
n\left(\max_{u\in V}y_u-\min_{u\in V}y_u\right),
\]
we have
\begin{equation}\label{eq:isotonic-as-tvd}
\operatorname{ISO}(y,D)
=
\operatorname{TVD}
\left(\widetilde y^{(\lambda)},\widetilde G,
\frac{\lambda}{2}\right).
\end{equation}
\end{proposition}

\begin{remark}[The one-dimensional chain]
For the chain graph \(1\to2\to\cdots\to n\), with \(n\geq2\), the
out-degree minus in-degree is \(1\) at the first vertex, \(-1\)
at the last vertex, and zero elsewhere. Thus,
\[
\widetilde y^{(\lambda)}
=
(y_1-\lambda/2,\;y_2,\ldots,y_{n-1},\;y_n+\lambda/2).
\]
Consequently, univariate isotonic regression can be computed by shifting
the first observation down by \(\lambda/2\), shifting the last
observation up by \(\lambda/2\), and applying univariate TVD with
tuning parameter \(\lambda/2\), provided
\(\lambda>n(\max_i y_i-\min_i y_i)\).
\end{remark}

As we will show in this article, applying our TVD formula to these transformed observations recovers
the classical isotonic representation
\eqref{eq:classical-isotonic-intro}. In this sense, our pointwise min--max representation for graph TVD is a strict and seemingly far-reaching generalization of Brunk's classical formula for
isotonic regression. We now present our general strategy for showing the reduction from isotonic regression to TVD as in Proposition~\ref{prop:isotonic-as-tvd}. It involves going through two intermediate classes of estimators, namely asymmetric TVD and nearly isotonic regression.

\subsubsection{Asymmetric TVD}

Let \(D=(V,A)\) be a finite DAG. For
\(\lambda_+,\lambda_-\geq0\), define the asymmetric TVD estimator by
\begin{equation}\label{eq:asymmetric-tvd}
\begin{aligned}
&\widehat\theta^{(\lambda_+,\lambda_-)}\\
&\quad=\argmin_{\theta\in\R^V}\Biggl\{
\frac12\sum_{u\in V}(y_u-\theta_u)^2\\
&\qquad\qquad+\lambda_+\sum_{(u,v)\in A}(\theta_v-\theta_u)_+\\
&\qquad\qquad+\lambda_-\sum_{(u,v)\in A}(\theta_v-\theta_u)_-
\Biggr\},
\end{aligned}
\end{equation}
where \(x_+=\max\{x,0\}\) and \(x_-=\max\{-x,0\}\).
Thus, increases and decreases along a directed edge receive
penalties \(\lambda_+\) and \(\lambda_-\), respectively.

The elementary identity
\[
\begin{gathered}
\lambda_+x_+ + \lambda_-x_- = a|x|+bx,\\
a=\frac{\lambda_++\lambda_-}{2},
\qquad
b=\frac{\lambda_+-\lambda_-}{2}
\end{gathered}
\]
allows us to express the asymmetric penalty as a symmetric TV
penalty plus a linear function of \(\theta\).
Indeed, let
\[
d_u=\deg_D^+(u)-\deg_D^-(u)
\]
denote the out-degree minus the in-degree of vertex \(u\).
Applying the identity to each edge and using
\[
\sum_{(u,v)\in A}(\theta_v-\theta_u)
=
-\sum_{u\in V}d_u\theta_u,
\]
we obtain
\begin{equation}\label{eq:asymmetric-penalty-decomposition}
\begin{aligned}
&\lambda_+\sum_{(u,v)\in A}(\theta_v-\theta_u)_+
+\lambda_-\sum_{(u,v)\in A}(\theta_v-\theta_u)_-\\
&\qquad=a\sum_{(u,v)\in A}|\theta_v-\theta_u|
-b\sum_{u\in V}d_u\theta_u.
\end{aligned}
\end{equation}
Absorbing the linear term into the squared loss by completing the square gives the following reduction.

\begin{lemma}[Asymmetric TVD as symmetric TVD]
\label{lem:asymmetric-to-symmetric}
Let \(\widetilde G\) be the underlying undirected graph of \(D\).
With \(a,b,d\) as above, define the shifted data vector
\[
\widetilde y_u=y_u+bd_u,
\qquad u\in V.
\]
Then
\begin{equation}\label{eq:asymmetric-as-symmetric}
\widehat\theta^{(\lambda_+,\lambda_-)}
=
\operatorname{TVD}(\widetilde y,\widetilde G,a).
\end{equation}
\end{lemma}

Thus, asymmetric TVD is symmetric TVD applied to suitably
transformed observations.

\subsubsection{Nearly Isotonic Regression}

Taking \(\lambda_+=0\) and \(\lambda_-=\lambda\) in
\eqref{eq:asymmetric-tvd} gives the nearly isotonic estimator
\begin{equation}\label{eq:nearly-isotonic}
\begin{aligned}
&\operatorname{NearIso}(y,D,\lambda)\\
&\quad=\argmin_{\theta\in\R^V}\Biggl\{
\frac12\sum_{u\in V}(y_u-\theta_u)^2\\
&\qquad\qquad+\lambda\sum_{(u,v)\in A}(\theta_v-\theta_u)_-
\Biggr\}.
\end{aligned}
\end{equation}
The penalty consists of only violations of the isotonic inequalities
\(\theta_u\leq\theta_v\), with larger values of \(\lambda\)
imposing a greater cost on such violations.

Nearly isotonic regression in the univariate setting was introduced by
\citet{tibshirani2011nearly}, who developed an algorithm for its
entire solution path. Extensions to directed graphs have been
studied by \citet{minami2020estimating,pastukhov2024fused,pastukhov2026generalized}.

For this choice of penalty,
\[
a=\frac{\lambda}{2},
\qquad
b=-\frac{\lambda}{2}.
\]
Lemma~\ref{lem:asymmetric-to-symmetric} therefore gives the
following reduction.

\begin{corollary}[Nearly isotonic regression as TVD]
\label{cor:neariso-as-tvd}
Let \(\widetilde G\) be the underlying undirected graph of \(D\),
and let \(d_u=\deg_D^+(u)-\deg_D^-(u)\). Then, for every
\(\lambda\geq0\),
\begin{equation}\label{eq:neariso-as-tvd}
\begin{aligned}
&\operatorname{NearIso}(y,D,\lambda)\\
&\quad=\operatorname{TVD}
\left(y-\frac{\lambda}{2}d,\widetilde G,\frac{\lambda}{2}\right).
\end{aligned}
\end{equation}
\end{corollary}

The next lemma shows that nearly isotonic regression agrees
exactly with isotonic regression once the penalty exceeds
an explicit finite threshold.

\begin{lemma}[Exact agreement with isotonic regression]
\label{lem:neariso-equals-iso}
Let \(D=(V,A)\) be a finite DAG with \(n=|V|\geq1\), and set
\[
\lambda_0
=
n\left(\max_{u\in V}y_u-\min_{u\in V}y_u\right).
\]
For every \(\lambda>\lambda_0\),
\begin{equation}\label{eq:neariso-equals-iso}
\operatorname{NearIso}(y,D,\lambda)
=
\operatorname{ISO}(y,D).
\end{equation}
\end{lemma}

\begin{remark}
As \(\lambda\to\infty\), violations of the order constraints
become increasingly costly, suggesting convergence to isotonic
regression. Lemma~\ref{lem:neariso-equals-iso} strengthens this
limiting statement to exact agreement at a finite penalty.
\end{remark}

\begin{remark}
The threshold \(\lambda_0\) is only a sufficient bound:
exact agreement may occur for smaller values of \(\lambda\),
and we do not claim that \(\lambda_0\) is the smallest threshold
for a given \(y\).
\end{remark}
Combining Corollary~\ref{cor:neariso-as-tvd} with
Lemma~\ref{lem:neariso-equals-iso} proves
Proposition~\ref{prop:isotonic-as-tvd}.

\subsubsection{From Nested Sets to Upper and Lower Sets}
We have thus traced a chain of special cases: isotonic regression
is a special case of nearly isotonic regression, which is a special
case of asymmetric TVD, which in turn reduces to symmetric TVD
after a shift of the observations. A single TVD minmax formula
therefore unifies all four estimators.

Our TVD minmax representation can be shown to recover the classical
isotonic formula. The mechanism behind this is worth highlighting
here. After substituting the transformed observations and the
penalty parameter into the min--max representation, the
correction terms involving \(C_{I,J}\), combined with the data
shift, simplify considerably. Since the penalty is sufficiently large, any outer set \(J\)
with an incoming edge is ruled out: the inner choice \(I=J\)
produces a large positive correction, making the inner maximum
too large to attain the outer minimum over \(J\).
Thus only lower sets \(J=L\) remain. For such an \(L\), any inner
set \(I\) with an edge directed from \(I\) to \(L\setminus I\)
receives a large negative correction and is excluded by the
inner maximization. The surviving inner sets are therefore
upper sets relative to \(L\), precisely the intersections
\(I=L\cap U\) with \(U\) an upper set. On these pairs the combined
correction from the data shift and TV penalty vanishes, leaving
only the local averages \(\bar y_{L\cap U}\) as in the classical
formula. The max--min representation follows by reversing the
roles of upper and lower sets. The full proof is given in Section~\ref{sec:proof-isotonic-as-tvd}.

\subsection{General convex losses and nonunique solutions}
\label{subsec:convex}

Our interest in extending the pointwise representation beyond squared
loss began with quantile regression. Given the importance of quantile
regression \citep{koenker1978regression} in statistics, it is natural to ask for an analogous
description of TVD for quantile regression. For \(0<\tau<1\), let
\(\rho_\tau(r)=\tau r_++(1-\tau)r_-\). In the univariate setting, the quantile TVD problem is defined as
\begin{equation}\label{eq:chain-quantile-tvd}
\cS^\tau_\lambda
=\argmin_{\theta\in\R^n}
\left\{\sum_{i=1}^n\rho_\tau(y_i-\theta_i)
+\lambda\sum_{i=1}^{n-1}|\theta_{i+1}-\theta_i|\right\}.
\end{equation}
Unlike squared loss, quantile loss is not strictly convex, so the
minimizer need not be unique. Recall that a function $f:C\to\mathbb R$ on a convex set
$C\subseteq\mathbb R^d$ is \emph{strictly convex} if
\[
f((1-t)x+ty)<(1-t)f(x)+tf(y)
\]
for every distinct $x,y\in C$ and every $t\in(0,1)$.
Strict convexity guarantees uniqueness of a minimizer whenever
one exists.

This raises a basic question: in the
absence of a unique fitted value, what could a pointwise minmax
representation characterize?

In the one dimensional case, \citet{ghoshal2026exact} answer this question by characterizing the set of possible fitted values at each coordinate. To state their
result in our notation, let \(\mathcal I_i\) be the intervals containing
\(i\), and compute \(C_{I,J}\) for the $1D$ chain graph as in \eqref{eq:boundary}.
For any nonempty set \(I\), write \(y_{I,(k)}\) for its \(k\)th order
statistic, with the conventions
\[
y_{I,(k)}=-\infty\quad(k\leq0),\qquad
y_{I,(k)}=+\infty\quad(k>|I|).
\]
Then,
\begin{equation}\label{eq:chain-quantile-interval}
\{\theta_i:\theta\in\cS^\tau_\lambda\}
=[\lo_i^{\,\tau},\hi_i^{\,\tau}],
\end{equation}
where
\begin{align}
\lo_i^{\,\tau}
&=\max_{J\in\mathcal I_i}
\min_{\substack{I\in\mathcal I_i\\I\subseteq J}}
y_{I,(\lceil\tau|I|+\lambda C_{I,J}\rceil)},
\label{eq:chain-quantile-lower}\\
\hi_i^{\,\tau}
&=\min_{J\in\mathcal I_i}
\max_{\substack{I\in\mathcal I_i\\I\subseteq J}}
y_{I,(\lfloor\tau|I|-\lambda C_{I,J}\rfloor+1)}.
\label{eq:chain-quantile-upper}
\end{align}
Thus, the possible fitted values at coordinate $i$ form an interval
whose left and right endpoints are, respectively, a maxmin and a
minmax of corrected local $\tau$-quantiles over the inner interval
$I$. The quantile ranks are adjusted by the same boundary correction
$C_{I,J}$ that appears in the squared-loss representation of
\citet{chatterjee2026minmax} and its graph extension above.

 Moreover, both endpoint vectors
\(\lo^{\,\tau}=(\lo_i^{\,\tau})_{i=1}^n\) and
\(\hi^{\,\tau}=(\hi_i^{\,\tau})_{i=1}^n\) belong to
\(\cS^\tau_\lambda\). Thus, the endpoints are attained simultaneously
at every coordinate by the componentwise smallest and largest
solutions, respectively.

This result suggests three further questions. Does the interval
description persist for an arbitrary convex loss? Can its endpoints
be expressed by a common minmax/maxmin formula that includes both
squared loss and quantile loss? And can such a characterization
hold on any finite graph? The answer to all these questions is yes.
The interval property follows from convexity and compactness; the
substantive issue is to identify its endpoints explicitly and show
that they are attained simultaneously by two extremal solutions.

Let \(G=(V,E)\) be a finite undirected graph, and let
\(\ell_u:\R\to\R\), \(u\in V\), be finite convex functions, with
dependence on \(y_u\) left implicit. For \(\lambda\geq0\), define the generalized graph TVD problem
\begin{equation}\label{eq:general-objective}
F_\lambda(\theta)
=\sum_{u\in V}\ell_u(\theta_u)+\lambda\TV_G(\theta),
\qquad
\cS_\lambda=\argmin_{\theta\in\R^V}F_\lambda(\theta).
\end{equation}
We assume that \(F_\lambda\) is coercive, meaning that
\(F_\lambda(\theta)\to+\infty\) as \(\|\theta\|\to\infty\).
This ensures that the solution set \(\cS_\lambda\) is nonempty
and compact. A simple sufficient condition is that each
\(\ell_u\) is coercive, that is,
\(\ell_u(x)\to+\infty\) as \(|x|\to\infty\): since
\(\TV_G\geq0\), we then have
\[
F_\lambda(\theta)\geq
\sum_{u\in V}\ell_u(\theta_u)\to+\infty
\qquad\text{as }\|\theta\|\to\infty.
\]
Individual coercivity is not necessary, however; see Remark~\ref{rem:logistic}. We emphasize
that, beyond finiteness and convexity of the losses, our only
requirement is coercivity of the full objective. This is a mild regularity condition, satisfied by squared,
quantile, and Huber losses, and by logistic TVD under the
conditions in Remark~\ref{rem:logistic}. In particular, neither
differentiability nor strict convexity is required.

To define the local quantities in our representation, we first
recall some facts from univariate convex analysis. Let
$h:\R\to\R$ be a finite convex function. Then, $h$ is continuous,
and its left and right derivatives, denoted by $h'_-$ and $h'_+$,
exist and are finite everywhere. Both derivatives are
nondecreasing and satisfy $h'_-(t)\leq h'_+(t)$ for every
$t\in\R$; see \citet[Theorem~24.1]{rockafellar1997convex}.
Define
\begin{equation}\label{eq:scalar-endpoints}
\begin{aligned}
r_-(h)&=\inf\{t:h'_+(t)\geq0\},\\
r_+(h)&=\sup\{t:h'_-(t)\leq0\},
\end{aligned}
\end{equation}
with \(\inf\varnothing=+\infty\) and
\(\sup\varnothing=-\infty\).

\begin{remark}\label{rem:convex_fn_min_endpoint}
    By convexity and continuity, the set of minimizers of $h$, when
nonempty, is a closed interval, possibly unbounded or a singleton.
If this interval is bounded, its left and right endpoints are
$r_-(h)$ and $r_+(h)$, respectively. See Section~\ref{subsec:proof-convex-endpoint} for a proof.
\end{remark}

\begin{remark}
The same endpoint interpretation holds when the minimizing interval
is unbounded, with endpoints understood in the extended real line.
For example, if $h(t)=t_+$, then $\argmin h=(-\infty,0]$, so
$r_-(h)=-\infty$ and $r_+(h)=0$.
The quantities $r_-(h)$ and $r_+(h)$ remain well defined even when
the infimum is not attained. For instance, if $h(t)=e^t$, there
is no minimizer, but $r_-(h)=r_+(h)=-\infty$, reflecting that the
infimum is approached as $t\to-\infty$. Similarly, $h(t)=e^{-t}$
gives $r_-(h)=r_+(h)=+\infty$.
\end{remark}
For nonempty nested sets \(I\subseteq J\), set
\begin{align}
q^-_{I,J}
&=r_-\!\left(t\mapsto
\sum_{u\in I}\ell_u(t)-\lambda C_{I,J}t\right),
\label{eq:lower-scalar}\\
q^+_{I,J}
&=r_+\!\left(t\mapsto
\sum_{u\in I}\ell_u(t)+\lambda C_{I,J}t\right).
\label{eq:upper-scalar}
\end{align}
The quantities \(q^-_{I,J}\) and \(q^+_{I,J}\) describe local fits (on I),
the fitting objective \(\sum_{u\in I}\ell_u(t)\) is modified
by a linear function with slope $-\lambda C_{I,J}$ or $+\lambda C_{I,J}$, where $C_{I,J}$ is the correction factor defined in~\eqref{eq:boundary}.

More precisely, \(q^-_{I,J}\) is the smallest minimizer of
\(\sum_{u\in I}\ell_u(t)-\lambda C_{I,J}t\), while
\(q^+_{I,J}\) is the largest minimizer of
\(\sum_{u\in I}\ell_u(t)+\lambda C_{I,J}t\).
If either extreme minimizer does not exist, the corresponding
quantity is interpreted using the extended-real definitions
of \(r_-\) and \(r_+\) above.

For example, under squared loss
\(\ell_u(t)=\frac12(y_u-t)^2\), these quantities reduce to
\[
q^-_{I,J}=\bar y_I+\frac{\lambda C_{I,J}}{|I|},
\qquad
q^+_{I,J}=\bar y_I-\frac{\lambda C_{I,J}}{|I|}.
\]
Thus, they generalize the corrected local averages appearing
in the squared-loss representation.

Similarly, under quantile loss 
\(\ell_u(t)=\rho_\tau(y_u-t)\), it can be checked that
\[
q^-_{I,J}
=y_{I,(\lceil\tau|I|+\lambda C_{I,J}\rceil)},
\qquad
q^+_{I,J}
=y_{I,(\lfloor\tau|I|-\lambda C_{I,J}\rfloor+1)},
\]
with the order-statistic conventions given above. Thus, the same
construction yields corrected local averages for squared loss
and corrected local quantiles for quantile loss. We now state our main result.

\begin{theorem}[Pointwise representation under convex loss]
\label{thm:convex}
Let \(G=(V,E)\) be a finite undirected graph with
\(V\ne\varnothing\), let \(\ell_u:\R\to\R\), \(u\in V\), be finite convex functions,
and let \(\lambda\geq0\). Suppose that \(F_\lambda\) in
\eqref{eq:general-objective} is coercive. Then, the following statements hold.
\begin{enumerate}
\item[(i)] For every \(i\in V\), the set of possible fitted values
at coordinate \(i\) is a nonempty compact interval:
\begin{equation}\label{eq:coordinate-interval}
\{\theta_i:\theta\in\cS_\lambda\}=[\lo_i,\hi_i].
\end{equation}

\item[(ii)] The endpoints are given by
\begin{align}
\lo_i
&=\max_{J\in\cC_i}
\min_{\substack{I\in\cC_i\\I\subseteq J}}q^-_{I,J},
\label{eq:convex-lower}\\
\hi_i
&=\min_{J\in\cC_i}
\max_{\substack{I\in\cC_i\\I\subseteq J}}q^+_{I,J}.
\label{eq:convex-upper}
\end{align}
Both formulas remain valid when \(I\) and \(J\) range over all
vertex sets satisfying \(i\in I\subseteq J\), without connectedness
restrictions. The final minmax/maxmin values are finite, although individual
\(q^-_{I,J}\) and \(q^+_{I,J}\) may be infinite.

\item[(iii)] The solution set \(\cS_\lambda\) is a lattice under
pointwise minimum and maximum, i.e., if \(\theta,\eta\in\cS_\lambda\),
then \(\theta\wedge\eta\) and \(\theta\vee\eta\) also belong to
\(\cS_\lambda\), where the two operations are taken coordinatewise.
Together with compactness, this implies that \(\cS_\lambda\) has a
coordinatewise smallest solution \(\lo=(\lo_i)_{i\in V}\) and a
coordinatewise largest solution \(\hi=(\hi_i)_{i\in V}\).
These solutions attain the lower and upper endpoints, respectively,
simultaneously at every coordinate, and every
\(\theta\in\cS_\lambda\) satisfies
\[
\lo_i\leq\theta_i\leq\hi_i\qquad\text{for all }i\in V.
\]

\item[(iv)] If \(F_\lambda\) is strictly convex, in particular if
every \(\ell_u\) is strictly convex, then \(\lo=\hi\), and
\(\cS_\lambda\) consists of a unique solution.
\end{enumerate}
\end{theorem}

\begin{remark}
Although every value in $[\lo_i,\hi_i]$ is attained at
coordinate $i$ by some solution, arbitrary choices from
these intervals need not combine to form a solution.
In particular, $\cS_\lambda$ need not equal the Cartesian
product $\prod_{i\in V}[\lo_i,\hi_i]$.
The theorem nevertheless guarantees that the vectors
$\lo$ and $\hi$ are both solutions.
\end{remark}

The above theorem gives a complete characterization of each coordinate projection of the solution
set of the generalized TVD problem on graphs.

For squared loss, part~\textup{(iv)} applies, and
Theorem~\ref{thm:convex} recovers Theorem~\ref{thm:squared}
as a special case. The generality of the result lies in its
applicability to arbitrary finite, separable convex losses:
neither strict convexity nor differentiability is required,
only the mild additional assumption that the full objective is coercive. 

\subsubsection{Intuition behind the general representation}
\label{subsec:convex-heuristic}

We briefly explain the reasoning that led us to
Theorem~\ref{thm:convex}. Our starting point was the case
of differentiable, strictly convex losses, where the minimizer is unique and the loss terms enter
the optimality conditions through their ordinary derivatives. Arguments similar to those for squared
loss suggest retaining the same max--min and min--max
structure, with the corrected local averages replaced by
solutions of the scalar equations
\[
\sum_{u\in I}\ell'_u(t)=\lambda C_{I,J}
\qquad\text{and}\qquad
\sum_{u\in I}\ell'_u(t)=-\lambda C_{I,J},
\]
respectively. When these solutions exist for every admissible
pair $I,J$, the fitted value is the max--min of the solutions
to the first equation, or equivalently the min--max of those
to the second. For squared loss, these solutions are precisely
the corrected averages in Theorem~\ref{thm:squared}.
For losses of the form $\ell_u(t)=b(t)-y_ut$, with $b$
differentiable and strictly convex, as arising in exponential
families, the equations reduce to
\[
b'(t)=\bar y_I\pm\frac{\lambda C_{I,J}}{|I|}.
\]
Since $b'$ is strictly increasing, this yields the same
max--min and min--max formulas as in the squared-loss case
for the transformed fitted values $b'(\widehat\theta_i)$;
see Remark~\ref{rem:other-losses}.

For general convex losses, the minimizer need no longer
be unique. Under our assumptions, however, the solution
set $\cS_\lambda$ is convex and compact, so the possible
fitted values at each coordinate form a nonempty compact
interval. We therefore sought representations of the
smallest and largest possible fitted values.

A key step was to recognize the scalar equations above as
the optimality conditions for minimizing
\[
f_{I,J}(t)
=\sum_{u\in I}\ell_u(t)-\lambda C_{I,J}t,
\qquad
\widetilde f_{I,J}(t)
=\sum_{u\in I}\ell_u(t)+\lambda C_{I,J}t.
\]
These local optimization problems remain meaningful without
differentiability or strict convexity. Their minimizers may
be nonunique, but whenever they exist, they form an interval.

To see how this suggests a representation of the extremal
fitted values, fix a vertex $i$ and consider the min--max
operator
\[
T_i(a)
=
\min_{J\in\cC_i}
\max_{\substack{I\in\cC_i\\I\subseteq J}}a_{I,J}.
\]
This operator is nondecreasing in each input: if
$a_{I,J}\leq b_{I,J}$ for every admissible pair $(I,J)$,
then $T_i(a)\leq T_i(b)$. The corresponding max--min
operator has the same monotonicity property. Thus, when each input ranges over a compact interval, choosing all their
lower endpoints gives the smallest output, while choosing
all their upper endpoints gives the largest.

This suggested that the smallest possible fitted value
should be the max--min of the smallest minimizers of
$f_{I,J}$, while the largest possible fitted value should
be the min--max of the largest minimizers of
$\widetilde f_{I,J}$. 
This was the intuition behind the definitions of
$q^-_{I,J}$ and $q^+_{I,J}$. Formulating these quantities
carefully also accommodates local problems with unbounded
solution sets or no minimizer. Theorem~\ref{thm:convex} confirms this intuition, extending the
familiar formulas for a unique solution to characterize the
smallest and largest possible fitted values.

We next apply this result to quantile TVD on graphs.

\subsubsection{Quantile TVD on graphs}

For \(0<\tau<1\), consider quantile TVD on a finite undirected
graph \(G=(V,E)\):
\begin{equation}\label{eq:graph-quantile-tvd}
\cS^\tau_\lambda
=\argmin_{\theta\in\R^V}
\left\{\sum_{u\in V}\rho_\tau(y_u-\theta_u)
+\lambda\TV_G(\theta)\right\}.
\end{equation}
The quantile losses are finite and convex, and the objective is
coercive. Applying Theorem~\ref{thm:convex} therefore gives the
following result.

\begin{proposition}[Pointwise representation for graph quantile TVD]
\label{prop:quantile}
Let \(G=(V,E)\) be a finite undirected graph with
\(V\ne\varnothing\). For every \(y\in\R^V\), \(0<\tau<1\),
and \(\lambda\geq0\), the following statements hold.
\begin{enumerate}
\item[(i)] For every \(i\in V\), the set of possible fitted values
at coordinate \(i\) is a nonempty compact interval:
\[
\{\theta_i:\theta\in\cS^\tau_\lambda\}
=[\lo_i^{\,\tau},\hi_i^{\,\tau}].
\]

\item[(ii)] The endpoints are given by
\begin{align}
\lo_i^{\,\tau}
&=\max_{J\in\cC_i}
\min_{\substack{I\in\cC_i\\I\subseteq J}}
y_{I,(\lceil\tau|I|+\lambda C_{I,J}\rceil)},
\label{eq:graph-quantile-lower}\\
\hi_i^{\,\tau}
&=\min_{J\in\cC_i}
\max_{\substack{I\in\cC_i\\I\subseteq J}}
y_{I,(\lfloor\tau|I|-\lambda C_{I,J}\rfloor+1)}.
\label{eq:graph-quantile-upper}
\end{align}
Both formulas remain valid when \(I\) and \(J\) range over all
vertex sets satisfying \(i\in I\subseteq J\), without connectedness
restrictions. The values are finite, although individual order
statistics in the formulas may be infinite under the conventions below.

\item[(iii)] The endpoint vectors
\(\lo^{\,\tau}=(\lo_i^{\,\tau})_{i\in V}\) and
\(\hi^{\,\tau}=(\hi_i^{\,\tau})_{i\in V}\) belong to
\(\cS^\tau_\lambda\). They are its coordinatewise smallest and
largest solutions, respectively, and satisfy
\[
\lo_i^{\,\tau}\leq\theta_i\leq\hi_i^{\,\tau}
\qquad\text{for every }\theta\in\cS^\tau_\lambda
\text{ and every }i\in V.
\]
Thus, the lower and upper endpoints are attained simultaneously
at every coordinate by \(\lo^{\,\tau}\) and \(\hi^{\,\tau}\),
respectively.
\end{enumerate}
\end{proposition}

This proposition generalizes the result of \citet{ghoshal2026exact}
from the one-dimensional chain graph to arbitrary finite graphs.
Intervals on the chain graph are replaced by connected vertex sets,
while the corrected ranks retain the same form. Recall here that
\(y_{I,(k)}\) denotes the \(k\)th order statistic of
\(\{y_u:u\in I\}\) for \(1\leq k\leq |I|\), with the conventions
\[
y_{I,(k)}=-\infty\quad(k\leq0),\qquad
y_{I,(k)}=+\infty\quad(k>|I|).
\]

\begin{remark}[Further examples and connections]
\label{rem:other-losses}
Theorem~\ref{thm:convex} can be applied to canonical
exponential-family losses
\[
\ell_u(t)=b(t)-y_ut,
\]
where \(b:\R\to\R\) is differentiable and strictly convex,
provided the objective \(F_\lambda\) is coercive. Writing
\(g=(b')^{-1}\), the local quantities become
\begin{align*}
q^-_{I,J}
&=g\!\left(\bar y_I+\frac{\lambda C_{I,J}}{|I|}\right),\\
q^+_{I,J}
&=g\!\left(\bar y_I-\frac{\lambda C_{I,J}}{|I|}\right),
\end{align*}
with extended-real interpretations inherited from the scalar
endpoint definitions in~\eqref{eq:scalar-endpoints}.
Since \(g\) is increasing, it commutes with the maxima and
minima in the representation. Consequently,
\[
b'(\widehat\theta)
=\operatorname{TVD}(y,G,\lambda),
\]
where \(b'\) acts coordinatewise and the right-hand side denotes
squared-loss graph TVD with the same penalty parameter.
This recovers the equivalence between TVD under exponential-family
loss and TVD under squared loss on arbitrary graphs established
by \citet{sadhanala2024exponential}, and yields explicit pointwise
representations for the corresponding natural-parameter estimates.
Examples include logistic TVD, with
\(b(t)=\log(1+e^t)\), and Poisson TVD, with \(b(t)=e^t\);
see also \citet{vandegeer2020logistic} for theoretical
properties of logistic TVD.

Beyond exponential-family likelihoods, the theorem can be applied to
Huber-loss TVD, considered computationally by
\citet[Section~5.1]{polson2016mixtures}, who do not provide a
pointwise characterization of the fitted values.
Since Huber loss is not strictly convex, the solution need not
be unique. Theorem~\ref{thm:convex} gives exact pointwise
minmax/maxmin representations for the coordinatewise smallest
and largest solutions on arbitrary graphs, with the local
quantities obtained by substituting the Huber loss
into~\eqref{eq:lower-scalar} and~\eqref{eq:upper-scalar}.
\end{remark}

\begin{remark}[Coercivity for logistic TVD]\label{rem:logistic}
Individual coercivity of the losses is not necessary for Theorem~\ref{thm:convex}
to hold. Let \(G\) be connected, let \(y_u\in\{0,1\}\), and consider
\[
F_\lambda(\theta)
=\sum_{u\in V}\bigl\{\log(1+e^{\theta_u})-y_u\theta_u\bigr\}
+\lambda\TV_G(\theta),
\qquad \lambda>0.
\]
Although none of the individual logistic losses is coercive,
\(F_\lambda\) is coercive as long as the observations include
both a zero and a one. Hence Theorem~\ref{thm:convex} applies.
Moreover, \(F_\lambda\) is strictly convex, so the solution is
unique. A proof is given in Section~\ref{subsec:proof-logistic}.
\end{remark}

\subsubsection{Isotonic regression under convex loss}
\label{subsec:isotonic}

The connection to isotonic regression also extends to general convex
losses. For the isotonic result below, we assume that each individual
loss is coercive. Let \(D=(V,A)\) be a finite DAG with \(V\ne\varnothing\),
and use the isotonic cone \(\mathcal K_D\) and the lower and upper
sets from Section~1.2. Consider the generalized isotonic regression
problem
\begin{equation}\label{eq:isotonic-convex}
\cS_{\mathrm{iso}}
=\argmin_{\theta\in\mathcal K_D}\sum_{u\in V}\ell_u(\theta_u),
\end{equation}
where the losses $\ell_u$ are finite convex functions that are coercive
on \(\mathbb R\). For nonempty \(B\subseteq V\), define
\begin{equation}\label{eq:isotonic-scalars}
q_B^-=r_-\!\left(t\mapsto\sum_{u\in B}\ell_u(t)\right),
\qquad
q_B^+=r_+\!\left(t\mapsto\sum_{u\in B}\ell_u(t)\right).
\end{equation}
These quantities give the smallest and largest optimal constant
fits (w.r.t the losses) on \(B\).

\begin{theorem}[Convex-loss isotonic representation]
\label{thm:isotonic}
Under the assumptions above, the following statements hold.
\begin{enumerate}
\item[(i)] For every \(i\in V\), the set of possible fitted values
at coordinate \(i\) is a nonempty compact interval:
\[
\{\theta_i:\theta\in\cS_{\mathrm{iso}}\}
=[\lo_i^{\,\mathrm{iso}},\hi_i^{\,\mathrm{iso}}].
\]

\item[(ii)] The endpoints are given by
\begin{align}
\lo_i^{\,\mathrm{iso}}
&=\max_{U\in\cU_i}\min_{L\in\cL_i}q^-_{L\cap U},
\label{eq:isotonic-lower}\\
\hi_i^{\,\mathrm{iso}}
&=\min_{L\in\cL_i}\max_{U\in\cU_i}q^+_{L\cap U}.
\label{eq:isotonic-upper}
\end{align}

\item[(iii)] The solution set \(\cS_{\mathrm{iso}}\) is a lattice
under pointwise minimum and maximum: if
\(\theta,\eta\in\cS_{\mathrm{iso}}\), then
\(\theta\wedge\eta\) and \(\theta\vee\eta\) also belong to
\(\cS_{\mathrm{iso}}\).
Together with compactness, this implies that \(\cS_{\mathrm{iso}}\)
has a coordinatewise smallest solution
\(\lo^{\,\mathrm{iso}}=(\lo_i^{\,\mathrm{iso}})_{i\in V}\) and a
coordinatewise largest solution
\(\hi^{\,\mathrm{iso}}=(\hi_i^{\,\mathrm{iso}})_{i\in V}\).
These solutions attain the lower and upper endpoints, respectively,
simultaneously at every coordinate, and every
\(\theta\in\cS_{\mathrm{iso}}\) satisfies
\[
\lo_i^{\,\mathrm{iso}}\leq\theta_i\leq\hi_i^{\,\mathrm{iso}}
\qquad\text{for all }i\in V.
\]

\item[(iv)] If \(\theta\mapsto\sum_{u\in V}\ell_u(\theta_u)\)
is strictly convex on \(\mathcal K_D\), in particular if every
\(\ell_u\) is strictly convex, then
\(\lo^{\,\mathrm{iso}}=\hi^{\,\mathrm{iso}}\), and
\(\cS_{\mathrm{iso}}\) consists of a unique solution.
\end{enumerate}
\end{theorem}
For squared loss, \(q_B^-=q_B^+=\bar y_B\), recovering
\eqref{eq:classical-isotonic-intro}. For quantile loss,
\[
q_B^-=y_{B,(\lceil\tau|B|\rceil)},\qquad
q_B^+=y_{B,(\lfloor\tau|B|\rfloor+1)},
\]
so the theorem gives the extremal isotonic quantile fits.

The proof follows the same sequence of ideas as in the proof of the square loss case as described in
Section~\ref{sec:proof-isotonic-as-tvd}. Basically the result follows by reducing graph isotonic regression to graph
TVD through asymmetric TVD and nearly isotonic regression, and
then applying Theorem~\ref{thm:convex}.

Minmax representations for isotonic mean and quantile regression
are classical. Extensions to more general losses also have
substantial precedents.
\citet{robertson1980algorithms} obtain such
representations on partial orders under a Cauchy mean-value condition;
\citet[Theorem~4.4]{mosching2020monotone} give extremal-solution
formulas for coercive losses on a chain graph, including quantile loss;
and \citet{jordan2022characterizing} develop minmax characterizations
for set-valued identifiable functionals on arbitrary partial orders.
Theorem~\ref{thm:isotonic} provides a general formulation of
the isotonic minmax representation on arbitrary finite DAGs:
it allows arbitrary finite, convex, coercive losses at the
vertices and characterizes both extremal solutions without
assuming differentiability or strict convexity. This formulation
may be of independent interest in shape-restricted inference.
The key point is that our graph TVD representation strictly
generalizes the isotonic minmax framework: the representation for
isotonic regression under general convex losses follows as a
special case.

\section{Proof of Theorem~\ref{thm:squared}}
\label{sec:proof-squared}

Write $\widehat\theta=\widehat\theta^{(\lambda)}$. We will first prove \eqref{eq:square-maxmin}--\eqref{eq:square-minmax} over arbitrary sets $J,I$ containing $i$ and then justify why the
connected set restriction holds.

\medskip
\noindent\textbf{Step 1: Optimality Conditions}

Define the TVD objective function
$$f(\theta) := 
\frac12\sum_{u\in V}(y_u-\theta_u)^2
+\lambda\sum_{\{u,v\}\in E}|\theta_u-\theta_v|.$$

Since $f$ is convex, KKT conditions imply that $\widehat\theta$ is a minimizer of $f$ if and only if $0 \in \partial(f)(\widehat\theta)$, where $\partial(f)(\widehat\theta)$ is the set of subgradients of $f$ at $\widehat\theta.$ Now, let us express $f$ as
\[f(\theta)=\frac12\sum_{u\in V}(y_u-\theta_u)^2
+\lambda\sum_{\{u,v\}\in E}h_{uv}(\theta).\]
Here $h_{uv}(\theta)=|\theta_u-\theta_v|$.
Since both the loss function and $h_{uv}$ are convex and finite, one can apply the sum rule for subdifferentials. To this end, let us calculate the subgradient of the TVD objective function with respect to the coordinate $\theta_u$.
For the squared-error term,
\[
\frac{\partial}{\partial \theta_u}
\left\{
\frac12\sum_{v\in V}(\theta_v-y_v)^2
\right\}
=
\theta_u-y_u.
\]
For every edge $\{u,v\}\in E$, consider quantities satisfying
\[
s_{uv}
\in
\begin{cases}
\{\lambda\}, & \widehat\theta_v>\widehat\theta_u,\\[3pt]
[-\lambda,\lambda], & \widehat\theta_v=\widehat\theta_u,\\[3pt]
\{-\lambda\}, & \widehat\theta_v<\widehat\theta_u.
\end{cases}
\]
Here the ordering of the pair $(u,v)$ is important: $s_{uv}$ denotes the
quantity associated with the edge $\{u,v\}$ when viewed from $u$ toward $v$.
For the reverse ordering, we set $s_{vu}=-s_{uv}.$
Then, the contribution of the total variation penalty to the
subgradient of $f$ with respect to $\theta_u$ is
precisely $-\sum_{v:\{u,v\}\in E} s_{uv}.$

Consequently, $\widehat\theta$ is a minimizer if and only if there exist such edge variables satisfying
\begin{equation}\label{eq:optimality}
\widehat\theta_u-y_u
=
\sum_{v:\{u,v\}\in E} s_{uv},
\qquad u\in V.
\end{equation}

\medskip
\noindent\textbf{Step 2: A Set Version of the Optimality Conditions.}

Fix any nonempty set $I\subseteq V$. Summing the optimality conditions~\eqref{eq:optimality} over all vertices in $I$, we obtain
\[
\begin{aligned}
\sum_{u\in I}(\widehat\theta_u-y_u)
&=\sum_{u\in I}\sum_{\substack{v\in V\\ \{u,v\}\in E}} s_{uv}\\
&=\sum_{u\in I}\sum_{\substack{v\in I\\ \{u,v\}\in E}} s_{uv}
+\sum_{u\in I}\sum_{\substack{v\in I^c\\ \{u,v\}\in E}} s_{uv}\\
&=\sum_{u\in I}\sum_{\substack{v\in I^c\\ \{u,v\}\in E}} s_{uv}.
\end{aligned}
\]

The last equality follows because,
if $\{u,v\}\in E$ with $u,v\in I$, then the edge appears once as $s_{uv}$
and once as $s_{vu}$, and by definition, $s_{vu}=-s_{uv}.$

Now let us introduce the following notation. For two disjoint sets $A,B\subseteq V$, define
\[
\mathcal F(A,B)
=
\sum_{\substack{u\in A,\ v\in B\\ \{u,v\}\in E}} s_{uv}.
\]
With this notation,
\[
\sum_{u\in I}(\widehat\theta_u-y_u)
=
\mathcal F(I,I^c).
\]

Defining $\overline{\widehat\theta}_I=|I|^{-1}\sum_{u\in I}\widehat\theta_u$ and dividing by $|I|$ gives
\begin{equation}\label{eq:set-optimality}
\overline{\widehat\theta}_I
=
\bar y_I+\frac{\mathcal F(I,I^c)}{|I|}.
\end{equation}

The above set identity for all sets $I$ also characterizes the optimal solution $\widehat\theta$ and is the set version of the optimality conditions~\eqref{eq:optimality}.

\medskip


\medskip
\noindent\textbf{Step 3: The Upper Bound.}

We will first prove the following minmax upper bound for a fixed vertex $i \in V.$ 
\begin{equation}\label{eq:upbd}
\widehat\theta_i^{(\lambda)}
\leq
\min_{\substack{J\subseteq V\\ i\in J}}
\max_{\substack{I\subseteq J\\ i\in I}}
\left\{
\bar y_I-\frac{\lambda C_{I,J}}{|I|}
\right\}.
\end{equation}
Note that there is no restriction above for $J,I$ to be connected sets.

Let $J\subseteq V$ be any set containing $i$.
Define
\[
I=\{u\in J:\widehat\theta_u\geq \widehat\theta_i\}.
\]
Clearly, $i\in I$.

Since $\widehat\theta_u\geq \widehat\theta_i$ for every $u\in I$, we have by~\eqref{eq:set-optimality},
\[
\widehat\theta_i
\leq
\overline{\widehat\theta}_I
=
\bar y_I+\frac{\mathcal F(I,I^c)}{|I|}.
\]

In view of the preceding display, to upper bound $\widehat\theta_i$ it remains to upper bound $\mathcal F(I,I^c)$ which is the sum of the quantities $s_{uv}$ over all edges from $I$ to $I^c$. Since each $s_{uv}\leq \lambda$, a trivial upper bound would be obtained by replacing every such term by $+\lambda$.

However, we can do better. The key observation is that, for certain boundary edges, the definition of $s_{uv}$ forces
$s_{uv}=-\lambda.$ Using the value $-\lambda$ on these edges gives an upper bound that is at least as sharp. Thus the main question is: for which edges $\{u,v\}$ with $u\in I$ and $v\in I^c$ can we determine that $s_{uv}=-\lambda$?

For edges with $u\in I$, $v\in J\setminus I$, and $\{u,v\}\in E$, the definition of $I$ gives

$$
\widehat\theta_u\geq \widehat\theta_i>\widehat\theta_v.
$$

Hence

$$
s_{uv}=-\lambda,
$$

and therefore

$$
\mathcal F(I,J\setminus I)
=
-\lambda |E(I,J\setminus I)|.
$$

On the other hand, for edges $\{u,v\}\in E$ with $u\in I$ and $v\in J^c$, we have no ordering information between $\widehat\theta_u$ and $\widehat\theta_v$. We therefore use the trivial bound $s_{uv}\leq\lambda$, which yields

$$
\mathcal F(I,J^c)
\leq
\lambda |E(I,J^c)|.
$$

Combining the two bounds gives
\[
\begin{aligned}
\mathcal F(I,I^c)&=\mathcal F(I,J\setminus I)+\mathcal F(I,J^{c})\\
&\leq-\lambda |E(I,J\setminus I)|
+\lambda |E(I,J^{c})|\\
&=-\lambda C_{I,J}.
\end{aligned}
\]
It follows that
\[
\widehat\theta_i
\leq
\overline{\widehat\theta}_I
=
\bar y_I+\frac{\mathcal F(I,I^c)}{|I|}
\leq
\bar y_I-\frac{\lambda C_{I,J}}{|I|}.
\]

Since the preceding argument holds for every set $J\ni i$, for each such $J$ there exists a set $I\subseteq J$ with $i\in I$ for which

$$
\widehat\theta_i
\leq
\bar y_I-\frac{\lambda C_{I,J}}{|I|}.
$$

Therefore,

$$
\widehat\theta_i
\leq
\max_{\substack{I\subseteq J\\ i\in I}}
\left\{
\bar y_I-\frac{\lambda C_{I,J}}{|I|}
\right\}
$$

for every $J\ni i$. Taking the minimum over all such $J$ proves~\eqref{eq:upbd}.

\medskip
\noindent\textbf{Step 4: The Lower Bound.}

We will now prove a matching lower bound for a fixed vertex $i\in V$:
\begin{equation}\label{eq:lowbd}
\widehat\theta_i^{(\lambda)}
\geq
\min_{\substack{J\subseteq V\\ i\in J}}
\max_{\substack{I\subseteq J\\ i\in I}}
\left\{
\bar y_I-\frac{\lambda C_{I,J}}{|I|}
\right\}.
\end{equation}
Again, there is no restriction above for $J,I$ to be connected sets.

Define
\[
J=\{u\in V:\widehat\theta_u\leq \widehat\theta_i\}.
\]
Clearly, $i\in J$.

Now let $I\subseteq J$ be any set containing $i$. Since
$\widehat\theta_u\leq \widehat\theta_i$ for every $u\in I$, we have by~\eqref{eq:set-optimality},
\[
\widehat\theta_i
\geq
\overline{\widehat\theta}_I
=
\bar y_I+\frac{\mathcal F(I,I^c)}{|I|}.
\]

In view of the preceding display, to lower bound $\widehat\theta_i$ it remains to now lower bound $\mathcal F(I,I^c).$ Since each $s_{uv}\geq -\lambda$, a trivial lower bound would be obtained by replacing every such term by $-\lambda$.

The key observation again is that, for certain edges, the definition of $s_{uv}$ forces
$s_{uv}=+\lambda$. Using the value $+\lambda$ on these edges gives a lower bound that is at least as sharp. 

We observe that for edges with $u\in I$, $v\in J^{c}$, and $\{u,v\}\in E$, the definition of $J$ gives
\[
\widehat\theta_u\leq \widehat\theta_i<\widehat\theta_v.
\]
Hence
\[
s_{uv}=+\lambda,
\]
and therefore
\[
\mathcal F(I,J^{c})
=
\lambda |E(I,J^{c})|.
\]

On the other hand, for edges $\{u,v\}\in E$ with $u\in I$ and $v\in J\setminus I$, we have no ordering information between $\widehat\theta_u$ and $\widehat\theta_v$. We therefore use the trivial bound $s_{uv}\geq -\lambda$, which yields
\[
\mathcal F(I,J\setminus I)
\geq
-\lambda |E(I,J\setminus I)|.
\]

Combining the two bounds gives
\[
\begin{aligned}
\mathcal F(I,I^c)&=\mathcal F(I,J\setminus I)+\mathcal F(I,J^{c})\\
&\geq-\lambda |E(I,J\setminus I)|
+\lambda |E(I,J^{c})|\\
&=-\lambda C_{I,J}.
\end{aligned}
\]
It follows that
\[
\widehat\theta_i
\geq
\overline{\widehat\theta}_I
=
\bar y_I+\frac{\mathcal F(I,I^c)}{|I|}
\geq
\bar y_I-\frac{\lambda C_{I,J}}{|I|}.
\]

Since the preceding argument holds for every set $I\subseteq J$ with $i\in I$, we obtain
\[
\widehat\theta_i
\geq
\max_{\substack{I\subseteq J\\ i\in I}}
\left\{
\bar y_I-\frac{\lambda C_{I,J}}{|I|}
\right\}.
\]
Since the particular set $J=\{u\in V:\widehat\theta_u\leq\widehat\theta_i\}$
is feasible in the sense that it contains $i$, it follows that
\[
\widehat\theta_i
\geq
\min_{\substack{J\subseteq V\\ i\in J}}
\max_{\substack{I\subseteq J\\ i\in I}}
\left\{
\bar y_I-\frac{\lambda C_{I,J}}{|I|}
\right\}.
\]

\medskip
\noindent\textbf{Step 5: The Connected Set Restriction.}

The previous two steps have shown that the same quantity is both an
upper and a lower bound for the fitted value $\widehat\theta_i$. Therefore,
we have established
\begin{equation}\label{eq:minmax-arbitrary}
\widehat\theta_i^{(\lambda)}
=
\min_{\substack{J\subseteq V\\ i\in J}}
\max_{\substack{I\subseteq J\\ i\in I}}
\left\{
\bar y_I-\frac{\lambda C_{I,J}}{|I|}
\right\}.
\end{equation}

This is a minmax representation over arbitrary subsets $I\subseteq J$
containing $i$. However, it turns out that both $I$ and $J$ can be
restricted to connected sets.

We first revisit the upper-bound argument in Step $3$. Let $J\ni i$ be any connected
set and define
\[
\widetilde I
=
J\cap\{u\in V:\widehat\theta_u\geq\widehat\theta_i\}.
\]
Now consider the subgraph induced by $\widetilde I$, and let $I$ be the
(uniquely defined) connected component of this subgraph containing $i$. Then $I$ is connected,
$i\in I$, and $I\subseteq J$.

We claim that for every edge $\{u,v\}\in E$ with
$u\in I$ and $v\in J\setminus I$, we have
\[
\widehat\theta_u\geq\widehat\theta_i>\widehat\theta_v.
\]
Indeed, if $v\in\widetilde I$, then the edge $\{u,v\}$ would imply that
$v$ belongs to the same connected component of the subgraph induced by
$\widetilde I$ as $u$, and hence $v\in I$, a contradiction. Therefore
$v\notin\widetilde I$, and so
\[
\widehat\theta_v<\widehat\theta_i.
\]
Thus
\[
s_{uv}=-\lambda
\]
for every edge from $I$ to $J\setminus I$, and the remainder of the
upper-bound argument proceeds exactly as before. Consequently,
\[
\widehat\theta_i
\leq
\min_{\substack{J\in\mathcal C\\ i\in J}}
\max_{\substack{I\in\mathcal C\\ I\subseteq J,\ i\in I}}
\left\{
\bar y_I-\frac{\lambda C_{I,J}}{|I|}
\right\}.
\]

For the reverse inequality, start as in the lower bound argument in Step $4$ with
\[
\widetilde J
=
\{u\in V:\widehat\theta_u\leq\widehat\theta_i\},
\]
and let $J$ be the connected component of the subgraph induced by
$\widetilde J$ containing $i$. In particular, $J$ is connected and
$i\in J$.

Now let $I\subseteq J$ be any connected set containing $i$. Consider an
edge $\{u,v\}\in E$ with $u\in I$ and $v\in J^{c}$. We must have
$v\notin\widetilde J$. Indeed, if $v\in\widetilde J$, then the edge
$\{u,v\}$ would connect $v$ to the connected component $J$ of
$\widetilde J$, implying $v\in J$, a contradiction. Hence
\[
\widehat\theta_v>\widehat\theta_i\geq\widehat\theta_u,
\]
and therefore
\[
s_{uv}=+\lambda.
\]

The remainder of the lower-bound argument proceeds exactly as before and yields
\[
\widehat\theta_i
\geq
\bar y_I-\frac{\lambda C_{I,J}}{|I|}
\]
for every connected $I\subseteq J$ containing $i$. Therefore,
\[
\widehat\theta_i
\geq
\min_{\substack{J\in\mathcal C\\ i\in J}}
\max_{\substack{I\in\mathcal C\\ I\subseteq J,\ i\in I}}
\left\{
\bar y_I-\frac{\lambda C_{I,J}}{|I|}
\right\}.
\]

Combining the two inequalities proves
\[
\widehat\theta_i^{(\lambda)}
=
\min_{\substack{J\in\mathcal C\\ i\in J}}
\max_{\substack{I\in\mathcal C\\ I\subseteq J,\ i\in I}}
\left\{
\bar y_I-\frac{\lambda C_{I,J}}{|I|}
\right\}.
\]

\medskip
\noindent\textbf{Step 6: The Maxmin Identity.}

The objective in \eqref{eq:graph-tvd} is invariant under the transformation
$(y,\theta)\mapsto(-y,-\theta)$. Since the minimizer is unique, it follows that
$
\widehat\theta^{(\lambda)}(-y)
=
-\widehat\theta^{(\lambda)}(y).
$
Applying the arbitrary-set and connected-set minmax identities just proved to the data vector $-y$ and then negating both sides exchanges the minimum and maximum, yielding the maxmin representation. 
\hfill$\square$

\section{Proof of Proposition~\ref{prop:isotonic-as-tvd}}
\label{sec:proof-isotonic-as-tvd}
The goal in this section is to illustrate how the classical isotonic regression minmax representations fall out from our TVD minmax representation. To this end, we prove Lemma~\ref{lem:asymmetric-to-symmetric},
Corollary~\ref{cor:neariso-as-tvd}, and
Lemma~\ref{lem:neariso-equals-iso} in order. The only nontrivial part here is Lemma~\ref{lem:neariso-equals-iso}.

\begin{proof}[Proof of Lemma~\ref{lem:asymmetric-to-symmetric}]
By \eqref{eq:asymmetric-penalty-decomposition}, the asymmetric
TVD objective equals
\[
\begin{aligned}
&\frac12\sum_{u\in V}(y_u-\theta_u)^2
+a\TV_{\widetilde G}(\theta)-b\sum_{u\in V}d_u\theta_u\\
&\qquad=
\frac12\sum_{u\in V}(y_u+bd_u-\theta_u)^2
+a\TV_{\widetilde G}(\theta)+c,
\end{aligned}
\]
where \(c=-\sum_{u\in V}(bd_uy_u+b^2d_u^2/2)\) is independent
of \(\theta\). Thus the two objectives have the same unique
minimizer, proving the lemma.
\end{proof}

\begin{proof}[Proof of Corollary~\ref{cor:neariso-as-tvd}]
Set \(\lambda_+=0\) and \(\lambda_-=\lambda\) in
Lemma~\ref{lem:asymmetric-to-symmetric}. Then \(a=\lambda/2\),
\(b=-\lambda/2\), and \(\widetilde y=y-\lambda d/2\), giving
\eqref{eq:neariso-as-tvd}.
\end{proof}

\begin{proof}[Proof of Lemma~\ref{lem:neariso-equals-iso}]
Write \(\widehat\theta^{[\lambda]}=\operatorname{NearIso}(y,D,\lambda)\)
and let
\[
y_{\min}=\min_{u\in V}y_u,\qquad
y_{\max}=\max_{u\in V}y_u.
\]
Fix \(\lambda>\lambda_0=n(y_{\max}-y_{\min})\) and \(i\in V\).

For disjoint sets \(S,T\subseteq V\), let \(A(S,T)\) denote the
directed edges from \(S\) to \(T\). Internal edges cancel when
we sum the (out minus in) degrees over \(I\), so
\[
\sum_{u\in I}d_u=|A(I,I^c)|-|A(I^c,I)|.
\]
Combining this with the definition of \(C_{I,J}\) on the
underlying undirected graph gives
\begin{equation}\label{eq:iso-proof-boundary}
\begin{aligned}
&C_{I,J}+\sum_{u\in I}d_u\\
&\qquad=2\bigl\{|A(I,J\setminus I)|-|A(J^c,I)|\bigr\}.
\end{aligned}
\end{equation}
Indeed, edges from \(I\) to \(J\setminus I\) contribute \(2\),
edges from \(J^c\) to \(I\) contribute \(-2\), and the other
boundary edges contribute zero.

Apply the arbitrary-set TVD formula \eqref{eq:minmax-arbitrary}
to the shifted observations \(y-\lambda d/2\) with tuning
parameter \(\lambda/2\). Corollary~\ref{cor:neariso-as-tvd}
and \eqref{eq:iso-proof-boundary} yield
\begin{equation}\label{eq:iso-proof-nested}
\begin{aligned}
\widehat\theta_i^{[\lambda]}
&=\min_{\substack{J\subseteq V\\i\in J}}
\max_{\substack{I\subseteq J\\i\in I}}
\Biggl\{\bar y_I\\
&\qquad-\lambda
\frac{|A(I,J\setminus I)|-|A(J^c,I)|}{|I|}\Biggr\}.
\end{aligned}
\end{equation}
We now simplify this expression using
\(y_{\min}\leq\bar y_I\leq y_{\max}\) for every nonempty \(I\).

\medskip
\noindent\textbf{Step 1: Restricting the outer minimum to lower sets.}

If \(J\ni i\) is not lower, then \(|A(J^c,J)|\geq1\).
The inner choice \(I=J\) gives
\[
\bar y_J+\lambda\frac{|A(J^c,J)|}{|J|}
\geq y_{\min}+\frac{\lambda}{n}>y_{\max}.
\]
In contrast, for \(J=V\), every term in the inner maximum is
at most \(\bar y_I\leq y_{\max}\). Thus no non-lower set can
attain the outer minimum.

We may therefore restrict \(J\) to \(L\in\cL_i\). Since no
edge enters \(L\) from outside, \(A(L^c,I)=\varnothing\)
for every \(I\subseteq L\), and hence
\begin{equation}\label{eq:iso-proof-lower}
\widehat\theta_i^{[\lambda]}
=\min_{L\in\cL_i}
\max_{\substack{I\subseteq L\\i\in I}}
\left\{\bar y_I-\lambda\frac{|A(I,L\setminus I)|}{|I|}\right\}.
\end{equation}

\medskip
\noindent\textbf{Step 2: Expressing the inner maximum using upper sets.}

Fix \(L\in\cL_i\). If \(A(I,L\setminus I)\ne\varnothing\), then
\[
\bar y_I-\lambda\frac{|A(I,L\setminus I)|}{|I|}
\leq y_{\max}-\frac{\lambda}{n}<y_{\min}.
\]
Such an \(I\) cannot attain the inner maximum, since \(I=L\)
gives \(\bar y_L\geq y_{\min}\). Only sets with
\(A(I,L\setminus I)=\varnothing\) remain, and their correction
vanishes.

We claim that these sets are precisely \(I=L\cap U\), with
\(U\in\cU_i\). Indeed, no edge can go from \(L\cap U\) to
\(L\setminus U\), since it would leave the upper set \(U\).
Conversely, if \(i\in I\subseteq L\) and
\(A(I,L\setminus I)=\varnothing\), set \(U=I\cup L^c\).
Then \(L\cap U=I\). An edge from \(U\) to \(U^c=L\setminus I\)
could start neither in \(I\), by assumption, nor in \(L^c\),
since \(L\) is lower. Thus \(U\) is upper and contains \(i\).

Substituting this characterization into
\eqref{eq:iso-proof-lower} gives
\begin{equation}\label{eq:iso-proof-minmax}
\widehat\theta_i^{[\lambda]}
=\min_{L\in\cL_i}\max_{U\in\cU_i}\bar y_{L\cap U},
\qquad \lambda>\lambda_0.
\end{equation}

\medskip
\noindent\textbf{Step 3: The maxmin identity.}

Let \(D^{\mathrm{rev}}\) be obtained by reversing every edge of
\(D\). The transformation
\((y,\theta,D)\mapsto(-y,-\theta,D^{\mathrm{rev}})\)
preserves the nearly isotonic objective. By uniqueness,
\[
\operatorname{NearIso}(-y,D^{\mathrm{rev}},\lambda)
=-\operatorname{NearIso}(y,D,\lambda).
\]
Reversal interchanges lower and upper sets, while negation
leaves the data range unchanged. Applying
\eqref{eq:iso-proof-minmax} to \(-y\) on \(D^{\mathrm{rev}}\)
and negating both sides therefore gives
\[
\widehat\theta_i^{[\lambda]}
=\max_{U\in\cU_i}\min_{L\in\cL_i}\bar y_{L\cap U}.
\]

\medskip
\noindent\textbf{Step 4: Exact agreement with isotonic regression.}

The right-hand side of \eqref{eq:iso-proof-minmax} does not
depend on \(\lambda\). Thus there is a fixed vector \(m\)
such that \(\widehat\theta^{[\lambda]}=m\) for every
\(\lambda>\lambda_0\). Write
\[
\begin{aligned}
f(\theta)&=\frac12\sum_{u\in V}(y_u-\theta_u)^2,\\
P(\theta)&=\sum_{(u,v)\in A}(\theta_v-\theta_u)_-.
\end{aligned}
\]
Since \(m\) minimizes \(f+\lambda P\) for every such \(\lambda\),
comparison with the zero vector gives
\[
f(m)+\lambda P(m)\leq f(0),\qquad \lambda>\lambda_0.
\]
Letting \(\lambda\to\infty\) shows that \(P(m)=0\), so \(m\)
is isotonic. For any isotonic vector \(\theta\), we also have
\(P(\theta)=0\), and optimality gives \(f(m)\leq f(\theta)\).
Hence \(m\) is the unique isotonic least-squares estimator.
This proves the lemma and, together with Steps 2 and 3,
recovers the classical representation
\eqref{eq:classical-isotonic-intro}.
\end{proof}

Combining Lemma~\ref{lem:neariso-equals-iso} with
Corollary~\ref{cor:neariso-as-tvd} proves
Proposition~\ref{prop:isotonic-as-tvd}.

\section{Proofs for general convex losses}
\label{sec:proof-convex}

We prove Theorem~\ref{thm:convex}, Proposition~\ref{prop:quantile},
and Theorem~\ref{thm:isotonic} in order.

\subsection{Proof of Theorem~\ref{thm:convex}}
\label{subsec:proof-general-tvd}

Write \(F=F_\lambda\) and \(\mathcal S=\cS_\lambda\).
We use the scalar endpoints \(r_\pm\) from
\eqref{eq:scalar-endpoints}. Let us first state the following elementary facts about convex maps that will be very useful in our proofs. 

\begin{lemma}[Scalar facts for finite convex functions]
\label{lem:convex-function-facts}\label{lem:endpoint-comparisons}
Let \(h:\R\to\R\) be a finite convex function and let \(t\in\R\).
\begin{enumerate}
\item If \(h'_-(t)\leq0\), then \(t\leq r_+(h)\).
\item If \(r_+(h)>t\), then \(h(t+\varepsilon)\leq h(t)\)
for all sufficiently small \(\varepsilon>0\).
\item For each \(\varepsilon>0\), the function
\(x\mapsto h(x+\varepsilon)-h(x)\) is nondecreasing.
\item If \(h'_+(t)<0\), then \(r_+(h)>t\).
\item If \(h'_-(t)>0\), then \(r_+(h)<t\).
\item If \(\widetilde h(x)=h(-x)\) for \(x\in\R\), then
\(r_+(\widetilde h)=-r_-(h)\).
\end{enumerate}
All relations involving \(r_+(h)\) and \(r_-(h)\) hold in the
extended real line: \(r_+(h)=+\infty\) is allowed in the first,
second, and fourth assertions, and \(r_+(h)=-\infty\) in the fifth.
\end{lemma}
The proof of Lemma~\ref{lem:convex-function-facts} is given in Section~\ref{sec:proof_aux}. Let us now get into the proof of Theorem~\ref{thm:convex}. For the ease of presentation, for any nonempty \(I\subseteq J\), let us define
\begin{equation}\label{eq:Phi-general}
\Phi_{I,J}(t):=\sum_{u\in I}\ell_u(t)+\lambda C_{I,J}t,
\end{equation}
and let
\begin{equation}\label{eq:Phi-minus-general}
\Phi^-_{I,J}(t):=\sum_{u\in I}\ell_u(t)-\lambda C_{I,J}t.
\end{equation}
Thus, \(q^+_{I,J}=r_+(\Phi_{I,J})\) and
\(q^-_{I,J}=r_-(\Phi^-_{I,J})\).

We first prove the interval assertion, then identify its endpoints.
We next show that the endpoint vectors are themselves solutions,
and finally establish uniqueness under strict convexity.

\medskip
\noindent\textbf{Step 1: Coordinate Intervals.}

Each finite convex loss $\ell_u$ is continuous on $\mathbb{R}$, so \(F\) is continuous on $\mathbb{R}^V$.
Coercivity of $F$ therefore makes the nonempty sublevel set
\[
K=\{\theta\in\R^V:F(\theta)\leq F(0)\}
\]
compact. The minimum of \(F\) on \(K\) is attained and is also
its global minimum. Hence, \(\mathcal S\) is nonempty and compact.
Being the set of minimizers of a convex map, \(\mathcal{S}\) is convex.

For each \(i\in V\), the coordinate projection of
\(\mathcal S\) is consequently a nonempty compact convex subset
of \(\R\), and thus a closed bounded interval. Define
\[
\lo_i=\min_{\theta\in\mathcal S}\theta_i,
\qquad
\hi_i=\max_{\theta\in\mathcal S}\theta_i.
\]
This proves \eqref{eq:coordinate-interval} and part~\textup{(i)}.
Next, we characterize these endpoints in terms of min--max and max--min formulae.

\medskip
\noindent\textbf{Step 2: Optimality Conditions.}

By the KKT conditions, a vector $\widehat\theta$ minimizes $F$ if and only if there
exist subgradients
\[
g_u\in\partial\ell_u(\widehat\theta_u)
\]
and edge variables $s_{uv}=-s_{vu}$ such that
\begin{equation}\label{eq:general-kkt}
g_u
=
\sum_{v:\{u,v\}\in E}s_{uv},
\qquad
s_{uv}\in
\begin{cases}
\{\lambda\},&\widehat\theta_v>\widehat\theta_u,\\
[-\lambda,\lambda],&\widehat\theta_v=\widehat\theta_u,\\
\{-\lambda\},&\widehat\theta_v<\widehat\theta_u.
\end{cases}
\end{equation}
Note that the dual variables $s_{uv}$ have the exact same definition as the duals in the squared loss case as they solely depend on the TV penalty. The difference from the squared loss case is that the sub-differential $\partial\ell_u(\widehat\theta_u)$ might not be singleton.
Indeed, since the losses are convex and continuous, the following interval characterization for the subgradient of $\ell_u$ holds (see proof of Theorem 24.1. of \cite{rockafellar1997convex}):
\begin{equation}\label{eqn:charac_subgradient}
    \partial\ell_u(\widehat\theta_u)=[ \ell'_{u,-}(\hat{\theta}_u),\, \ell'_{u,+}(\hat{\theta}_u)].
\end{equation}
For the squared loss case, both left and right derivative coincide and hence the sub-differential becomes a singleton.

For any nonempty $I\subseteq V$, summing these conditions
over $u\in I$ gives
\begin{align*}
\sum_{u\in I}g_u
&=
\sum_{u\in I}\sum_{v:\{u,v\}\in E}s_{uv}\\
&=
\sum_{u\in I}
\sum_{\substack{v\in I\\\{u,v\}\in E}}s_{uv}
+
\sum_{u\in I}
\sum_{\substack{v\in I^c\\\{u,v\}\in E}}s_{uv}\\
&=
\sum_{u\in I}
\sum_{\substack{v\in I^c\\\{u,v\}\in E}}s_{uv}.
\end{align*}
The internal edges cancel because $s_{uv}=-s_{vu}$.
Thus, with the boundary-flow notation from Section~\ref{sec:proof-squared},
\begin{equation}\label{eq:general-set-optimality}
\sum_{u\in I}g_u=\mathcal F(I,I^c).
\end{equation}

\medskip
\noindent\textbf{Step 3: The Upper Bound.}

Let $\widehat\theta$ be any solution and fix $i\in V$.
For any $J\subseteq V$ containing $i$, take
\[
I=J\cap\{u\in V:\widehat\theta_u\geq\widehat\theta_i\}.
\]
For every $u\in I$, it follows from the interval characterization of subdifferential in \eqref{eqn:charac_subgradient} and the monotonicity of the
left derivative of a finite convex function that
\[
g_u
\geq\ell'_{u,-}(\widehat\theta_u)
\geq\ell'_{u,-}(\widehat\theta_i),
\]
Consequently,
\[
\sum_{u\in I}\ell'_{u,-}(\widehat\theta_i)
\leq
\sum_{u\in I}g_u
=
\mathcal F(I,I^c).
\]

As in the squared-loss proof, every edge from $I$ to
$J\setminus I$ leads to a strictly smaller fitted value by definition of $I$.
Thus, its edge variable is $-\lambda$. On edges from $I$
to $J^c$, we use the trivial bound $s_{uv}\leq\lambda$. Therefore, it follows that
\[
\mathcal F(I,I^c)
\leq
-\lambda|E(I,J\setminus I)|
+\lambda|E(I,J^c)|
=
-\lambda C_{I,J}.
\]
Therefore,
\[
\sum_{u\in I}\ell'_{u,-}(\widehat\theta_i)
+\lambda C_{I,J}
\leq0,
\]
or equivalently,
\[
\Phi'_{I,J,-}(\widehat\theta_i)\leq0.
\]
By the definition of $r_+$,
\[
\widehat\theta_i\leq r_+(\Phi_{I,J})\leq \max_{\substack{I\subseteq J\\i\in I}}
r_+(\Phi_{I,J}).
\]

For every outer set $J\ni i$, we have found an inner set
$I\subseteq J$ containing $i$ such that the above inequality holds. Hence, taking minimum over the outer set $J$, we obtain
\begin{equation}\label{eq:general-upper-bound}
\widehat\theta_i
\leq
\min_{\substack{J\subseteq V\\i\in J}}
\max_{\substack{I\subseteq J\\i\in I}}
r_+(\Phi_{I,J}).
\end{equation}

\medskip
\noindent\textbf{Step 4: A Matching Bound for the Upper Endpoint.}

The preceding upper bound applies to every solution. To obtain a
matching bound for the upper endpoint, we now use its maximality
through a one-sided perturbation.

By compactness, choose a solution \(\widehat\theta\in\mathcal S\)
such that \(\widehat\theta_i=\hi_i\), and define
\[
J=\{u\in V:\widehat\theta_u\leq\widehat\theta_i\}.
\]
Take any $I\subseteq J$ containing $i$. We will show that
\[
\widehat\theta_i\geq r_+(\Phi_{I,J}).
\]

For $\varepsilon>0$, define the local perturbation
\[
\widetilde\theta_u
=
\begin{cases}
\widehat\theta_u+\varepsilon,&u\in I,\\
\widehat\theta_u,&u\notin I.
\end{cases}
\]
Since \(\widetilde\theta_i=\hi_i+\varepsilon>\hi_i\), by definition of $\hi_i$ the
perturbed vector cannot belong to \(\mathcal S\). Therefore, the perturbed vector should produce strictly larger objective, i.e., 
\begin{equation}\label{eq:strict-perturbation}
0<F(\widetilde\theta)-F(\widehat\theta).
\end{equation}

We now bound this change in the objective from above. For the loss term,
\begin{align*}
\sum_{u\in V}
\{\ell_u(\widetilde\theta_u)-\ell_u(\widehat\theta_u)\}
&=
\sum_{u\in I}
\{\ell_u(\widehat\theta_u+\varepsilon)
  -\ell_u(\widehat\theta_u)\}\\
&\leq
\sum_{u\in I}
\{\ell_u(\widehat\theta_i+\varepsilon)
  -\ell_u(\widehat\theta_i)\}.
\end{align*}
The last inequality uses $\widehat\theta_u\leq\widehat\theta_i$
on $I$ and the monotonicity of the discrete derivative
from Lemma~\ref{lem:convex-function-facts}.

Let us now track the amount of change in the penalty term after perturbation. Write
\[
\operatorname{TV}(\theta)
=
\sum_{\{u,v\}\in E}|\theta_u-\theta_v|.
\]
Edges internal to $I$ or to $I^c$ contribute no change.
Each edge from $I$ to $J\setminus I$ contributes at most
$\varepsilon$. Now, by definition of $J$, for an edge from $u\in I$ to $v\in J^c$, we have
\[
\widehat\theta_u\leq\widehat\theta_i<\widehat\theta_v.
\]
Thus, for all sufficiently small $\varepsilon>0$, the
contribution of this edge towards the change in the penalty is exactly $-\varepsilon$. Consequently,
\begin{align*}
\operatorname{TV}(\widetilde\theta)
-\operatorname{TV}(\widehat\theta)
&\leq
-\varepsilon|E(I,J^c)|
+\varepsilon|E(I,J\setminus I)|\\
&=\varepsilon C_{I,J}.
\end{align*}

Combining these bounds gives the key inequality
\begin{align}
F(\widetilde\theta)-F(\widehat\theta)
&\leq
\sum_{u\in I}
\{\ell_u(\widehat\theta_i+\varepsilon)
  -\ell_u(\widehat\theta_i)\}
+\lambda\varepsilon C_{I,J}\notag\\
&=
\Phi_{I,J}(\widehat\theta_i+\varepsilon)
-\Phi_{I,J}(\widehat\theta_i).
\label{eq:key-perturbation-bound}
\end{align}
Together with \eqref{eq:strict-perturbation}, this yields
\[
0<
\Phi_{I,J}(\widehat\theta_i+\varepsilon)
-\Phi_{I,J}(\widehat\theta_i)
\]
for every sufficiently small $\varepsilon>0$.
By the contrapositive statement of the second assertion of Lemma~\ref{lem:convex-function-facts}, this implies
\[
r_+(\Phi_{I,J})\leq\widehat\theta_i.
\]

Since this holds for every $I\subseteq J$ containing $i$, we conclude that
\[
\max_{\substack{I\subseteq J\\i\in I}}
r_+(\Phi_{I,J})
\leq\widehat\theta_i.
\]
Using this particular $J$ in the outer minimum gives
\[
\min_{\substack{J\subseteq V\\i\in J}}
\max_{\substack{I\subseteq J\\i\in I}}
r_+(\Phi_{I,J})
\leq\widehat\theta_i.
\]
Combining this with the upper bound proves
\[
\overline{\widehat\theta}_i
=
\min_{\substack{J\subseteq V\\i\in J}}
\max_{\substack{I\subseteq J\\i\in I}}
r_+(\Phi_{I,J}).
\]

\medskip
\noindent\textbf{Step 5: The Connected Set Restriction.}

The connected set restriction follows as in Section~\ref{sec:proof-squared}.

For the upper bound, let $J$ be connected and replace
\[
J\cap\{u:\widehat\theta_u\geq\widehat\theta_i\}
\]
by its connected component containing $i$.
Every edge from this component to $J\setminus I$ still leads
to a fitted value strictly below $\widehat\theta_i$.
The upper-bound argument therefore remains valid.

For the lower bound, take $J$ to be the connected component
containing $i$ of
\[
\{u:\widehat\theta_u\leq\widehat\theta_i\}.
\]
Every neighbor of $J$ outside $J$ has fitted value strictly
above $\widehat\theta_i$; otherwise it would belong to the
same component. The perturbation argument therefore applies
to every connected $I\subseteq J$ containing $i$.
This proves \eqref{eq:convex-upper}.

\medskip
\noindent\textbf{Step 6: The Maxmin Representation.}

Apply the minmax representation to the reflected losses
\[
\widetilde\ell_u(x)=\ell_u(-x).
\]
The reflected objective is \(\widetilde F(\theta)=F(-\theta)\),
which is again convex and coercive. Its solution set is
\(-\mathcal S\), so its upper endpoint at coordinate \(i\) is
\[
\max_{\eta\in-\mathcal S}\eta_i
=-\min_{\theta\in\mathcal S}\theta_i=-\lo_i.
\]
Moreover, its corresponding scalar function satisfies
\[
\widetilde\Phi_{I,J}(x)
=
\sum_{u\in I}\ell_u(-x)+\lambda C_{I,J}x
=
\Phi^-_{I,J}(-x).
\]
Since
\[
r_+\bigl(f(-\,\cdot)\bigr)=-r_-(f),
\]
negating the reflected minmax identity exchanges minimum
and maximum and gives \eqref{eq:convex-lower}.
The same argument applies with or without the connected
set restrictions. This completes the proof of part~\textup{(ii)}.

\medskip
\noindent\textbf{Step 7: Lattice Structure and Simultaneous Attainment.}

We first show that \(\mathcal S\) is a lattice under pointwise
minimum and maximum. For \(\theta,\eta\in\R^V\), write
\[
(\theta\wedge\eta)_u=\min\{\theta_u,\eta_u\},
\qquad
(\theta\vee\eta)_u=\max\{\theta_u,\eta_u\}.
\]
We use the elementary inequality
\[
\begin{multlined}
|\min(a,c)-\min(b,d)|+|\max(a,c)-\max(b,d)|\\
\leq |a-b|+|c-d|.
\end{multlined}
\]
To see this, note that if \(a\le c\) and \(b\le d\), or
\(a\geq c\) and \(b\geq d\), the left-hand side equals
\(|a-b|+|c-d|\). Otherwise, say \(a\le c\) and \(b>d\); then
the left-hand side is \(|a-d|+|c-b|\), and the claim follows
since \(a-d\) and \(c-b\) both lie between \(a-b\) and \(c-d\).
Applying this inequality to each edge and using separability
of the losses gives
\[
F(\theta\wedge\eta)+F(\theta\vee\eta)
\leq F(\theta)+F(\eta).
\]
If \(\theta,\eta\in\mathcal S\), both terms on the left are
at least the minimum of \(F\), so both
\(\theta\wedge\eta\) and \(\theta\vee\eta\) also belong to
\(\mathcal S\). This proves the lattice assertion.

By compactness, for each \(i\in V\) we can choose solutions
\(\theta^{i,-}\) and
\(\theta^{i,+}\) satisfying
\(\theta^{i,-}_i=\lo_i\) and \(\theta^{i,+}_i=\hi_i\).
Since \(V\) is finite, closure under coordinatewise minimum and
maximum implies that
\[
\bigwedge_{i\in V}\theta^{i,-}\in\mathcal S,
\qquad
\bigvee_{i\in V}\theta^{i,+}\in\mathcal S.
\]
At every coordinate \(u\), all solutions have values in
\([\lo_u,\hi_u]\), while \(\theta^{u,-}\) and \(\theta^{u,+}\)
attain the respective endpoints in $\mathcal S$. Thus, these two vectors are
exactly \(\lo\) and \(\hi\). This proves part~\textup{(iii)}.

\medskip
\noindent\textbf{Step 8: Uniqueness under Strict Convexity.}

If \(F\) is strictly convex, two distinct minimizers cannot exist, as
the midpoint of these two points would have objective value strictly below the
minimum. Since a minimizer exists, it is unique, and hence it follows that
\(\lo=\hi\). Finally, if every \(\ell_u\) is strictly convex,
then \(\theta\mapsto\sum_{u\in V}\ell_u(\theta_u)\) is strictly
convex. Adding the convex TV penalty preserves strict convexity,
which proves part~\textup{(iv)}.
\hfill$\square$

\subsection{Proof of Proposition~\ref{prop:quantile}}
\label{subsec:proof-quantile}

For \(\ell_u(t)=\rho_\tau(y_u-t)\), the one-sided derivatives are
\[
\ell'_{u,-}(t)=\mathbf1\{y_u<t\}-\tau,\qquad
\ell'_{u,+}(t)=\mathbf1\{y_u\leq t\}-\tau.
\]
Consequently, with \(m=|I|\),
\begin{align*}
q^-_{I,J}
&=\inf\left\{t:\#\{u\in I:y_u\leq t\}
                 \geq m\tau+\lambda C_{I,J}\right\}\\
&=y_{I,(\lceil m\tau+\lambda C_{I,J}\rceil)},\\
q^+_{I,J}
&=\sup\left\{t:\#\{u\in I:y_u<t\}
                 \leq m\tau-\lambda C_{I,J}\right\}\\
&=y_{I,(\lfloor m\tau-\lambda C_{I,J}\rfloor+1)}.
\end{align*}
These identities include ties and ranks outside \(1,\ldots,m\)
under the stated infinite-rank conventions. The ceiling and
floor-plus-one arise from the weak and strict inequalities,
respectively. Substitution into Theorem~\ref{thm:convex} proves
all assertions of the proposition.
\hfill\(\square\)

\subsection{Proof of Theorem~\ref{thm:isotonic}}
\label{subsec:proof-isotonic-convex}

Throughout this subsection, for \(\theta\in\R^V\) we write
\[
\begin{aligned}
f(\theta)&=\sum_{u\in V}\ell_u(\theta_u),\\
P(\theta)&=\sum_{(u,v)\in A}(\theta_v-\theta_u)_-,
\end{aligned}
\]
as in Step~4 of the proof of Lemma~\ref{lem:neariso-equals-iso},
with the squared loss replaced by a general separable convex loss.
Every summand of \(P\) is nonnegative, and the summand of an edge
\((u,v)\in A\) vanishes if and only if \(\theta_u\leq\theta_v\). Hence
\begin{equation}\label{eq:P-zero-iff-isotonic}
P\geq0\ \text{on }\R^V,
\qquad
P(\theta)=0\iff\theta\in\mathcal K_D.
\end{equation}
The standing assumptions of Theorem~\ref{thm:isotonic} imply that
every \(\ell_u\) is finite, convex, and coercive, and hence, 
\begin{equation}\label{eq:coercive-on-cone}
f(\theta^{(k)})\to+\infty
\ \text{ whenever }
\|\theta^{(k)}\|\to\infty, \;\;\theta^{(k)}\in \mathbb{R}^V.
\end{equation}
For \(\lambda\geq0\), the nearly isotonic objective under the
losses \(\ell_u\) is
\begin{equation}\label{eq:general-neariso-objective}
H_\lambda(\theta)=f(\theta)+\lambda P(\theta),
\qquad\theta\in\R^V,
\end{equation}
which reduces to the objective in \eqref{eq:nearly-isotonic}
when \(\ell_u(x)=\frac12(y_u-x)^2\).

The proof follows the squared-loss argument of
Section~\ref{sec:proof-isotonic-as-tvd}. There, isotonic regression
was reduced to TVD in three stages:
Lemma~\ref{lem:asymmetric-to-symmetric} and
Corollary~\ref{cor:neariso-as-tvd} expressed nearly isotonic
regression as TVD, Lemma~\ref{lem:neariso-equals-iso} identified
nearly isotonic regression with isotonic regression at a large
finite penalty, and the proof of Lemma~\ref{lem:neariso-equals-iso}
reduced the nested-set formula to lower and upper sets.
Lemma~\ref{lem:general-neariso-as-tvd} below is the counterpart of
the first stage, Lemma~\ref{lem:general-neariso-equals-iso} is the
counterpart of the second, and Steps~1--3 in the proof of
Theorem~\ref{thm:isotonic} are the counterparts of Steps~1--3 in
the proof of Lemma~\ref{lem:neariso-equals-iso}.

The next lemma is the convex-loss counterpart of
Lemma~\ref{lem:asymmetric-to-symmetric} and
Corollary~\ref{cor:neariso-as-tvd}. Under squared loss, the linear
term produced by the asymmetric penalty was absorbed into the data
by completing the square; for a general loss it is absorbed into
the loss itself.

\begin{lemma}[Nearly isotonic regression as generalized TVD]
\label{lem:general-neariso-as-tvd}
Let \(\widetilde G\) be the underlying undirected graph of \(D\),
let \(d_u=\deg_D^+(u)-\deg_D^-(u)\), and for \(\lambda\geq0\)
define
\begin{equation}\label{eq:tilted-losses}
\ell^{(\lambda)}_u(x)=\ell_u(x)+\frac{\lambda}{2}\,d_u\,x,
\qquad u\in V,\ x\in\R.
\end{equation}
Then each \(\ell^{(\lambda)}_u\) is finite and convex, and for every
\(\theta\in\R^V\),
\begin{equation}\label{eq:general-neariso-reduction}
H_\lambda(\theta)
=\sum_{u\in V}\ell^{(\lambda)}_u(\theta_u)
+\frac{\lambda}{2}\TV_{\widetilde G}(\theta).
\end{equation}
In other words, \(H_\lambda\) is the generalized graph TVD objective
\eqref{eq:general-objective} on the graph \(\widetilde G\), with
losses \(\ell^{(\lambda)}_u\) and penalty parameter \(\lambda/2\).
\end{lemma}

\begin{proof}
Each \(\ell^{(\lambda)}_u\) is the sum of a finite convex function
and a linear function, so it is finite and convex. Apply
\eqref{eq:asymmetric-penalty-decomposition} with \(\lambda_+=0\)
and \(\lambda_-=\lambda\), so that \(a=\lambda/2\) and
\(b=-\lambda/2\):
\[
\lambda P(\theta)
=\frac{\lambda}{2}\sum_{(u,v)\in A}|\theta_v-\theta_u|
+\frac{\lambda}{2}\sum_{u\in V}d_u\theta_u.
\]
Since \(D\) is acyclic, no two vertices are joined by directed edges
in both directions. The directed edges of \(D\) therefore correspond
one-to-one to the edges of \(\widetilde G\), and the first sum equals
\(\TV_{\widetilde G}(\theta)\). Adding \(f(\theta)\) to both sides
and grouping the linear term with the losses gives
\eqref{eq:general-neariso-reduction}.
\end{proof}

The next lemma is the convex-loss counterpart of
Lemma~\ref{lem:neariso-equals-iso}.

\begin{lemma}[Exact agreement with isotonic regression]
\label{lem:general-neariso-equals-iso}
Under the assumptions of Theorem~\ref{thm:isotonic}, the following
statements hold.
\begin{enumerate}
\item[(a)] The set \(\cS_{\mathrm{iso}}\) is nonempty, compact,
and convex. We write
\(f^*=\min_{\theta\in\mathcal K_D}f(\theta)\).
\item[(b)] There is a finite constant \(\gamma_*\geq0\) such that,
for every \(\theta\in\R^V\) and every \(\lambda\geq0\),
\begin{equation}\label{eq:exact-penalty-gap}
H_\lambda(\theta)-f^*\geq(\lambda-\gamma_*)P(\theta).
\end{equation}
\item[(c)] For every \(\lambda>\gamma_*\),
\begin{equation}\label{eq:general-exact-penalty-proof}
\argmin_{\theta\in\R^V}H_\lambda(\theta)=\cS_{\mathrm{iso}}.
\end{equation}
\end{enumerate}
\end{lemma}

\begin{proof}
\emph{Part (a).} Each \(\ell_u\) is finite and convex, hence
continuous, so \(f\) is continuous and convex on \(\R^V\). The cone
\(\mathcal K_D\) is closed and convex and contains \(0\). Consider
\[
K=\{\theta\in\mathcal K_D:f(\theta)\leq f(0)\}.
\]
This set is nonempty and closed. It is also bounded: otherwise it
would contain a sequence \(\theta^{(k)}\) with
\(\|\theta^{(k)}\|\to\infty\) and \(f(\theta^{(k)})\leq f(0)\),
contradicting \eqref{eq:coercive-on-cone}. Thus \(K\) is compact,
and \(f\) attains its minimum over \(K\). This minimum is the
minimum \(f^*\) of \(f\) over \(\mathcal K_D\), because
\(f>f(0)\geq f^*\) on \(\mathcal K_D\setminus K\). Hence,
\[
\cS_{\mathrm{iso}}=\{\theta\in\mathcal K_D:f(\theta)=f^*\}
\]
is nonempty. It is a closed subset of \(K\) due to continuity of $f$, hence it is compact, and it
is convex because both \(\mathcal K_D\) and \(f\) are convex.

\emph{Part (b).} Fix \(\theta^*\in\cS_{\mathrm{iso}}\). By
definition, \(\theta^*\) minimizes the finite convex function \(f\)
subject to the affine inequality constraints
\[
c_{uv}(\theta):=\theta_u-\theta_v\leq0,\qquad(u,v)\in A.
\]
These constraints are strictly feasible. Indeed, since \(D\) is
acyclic, it admits a topological ordering, that is, a bijection
\(\pi:V\to\{1,\ldots,n\}\) with \(\pi(u)<\pi(v)\) for every
\((u,v)\in A\), and the vector \(\theta^0=(\pi(u))_{u\in V}\)
satisfies \(c_{uv}(\theta^0)<0\) for every \((u,v)\in A\). The
optimal value \(f^*\) is finite. The Kuhn--Tucker theorem
\citep[Theorems~28.2 and~28.3]{rockafellar1997convex} therefore
provides multipliers \(\gamma_{uv}\geq0\), \((u,v)\in A\), such that
\[
\begin{gathered}
0\in\partial f(\theta^*)
+\sum_{(u,v)\in A}\gamma_{uv}\,\nabla c_{uv},\\
\gamma_{uv}\,c_{uv}(\theta^*)=0,\quad(u,v)\in A.
\end{gathered}
\]
Here \(\nabla c_{uv}=e_u-e_v\), where \(e_u\) is the \(u\)-th
standard basis vector of \(\R^V\). Moreover, since \(f\) is
separable, a vector \(g\) belongs to \(\partial f(\theta^*)\) if and
only if \(g_u\in\partial\ell_u(\theta^*_u)\) for every \(u\in V\). Reading the first condition coordinate by coordinate, we obtain for each \(u\in V\), there exists
\(g_u\in\partial\ell_u(\theta^*_u)\), such that
\begin{equation}\label{eq:general-iso-kkt-proof}
\begin{gathered}
g_u+\sum_{v:(u,v)\in A}\gamma_{uv}
-\sum_{v:(v,u)\in A}\gamma_{vu}=0,\quad u\in V,\\
\gamma_{uv}(\theta^*_u-\theta^*_v)=0,\quad(u,v)\in A.
\end{gathered}
\end{equation}
If \(A=\varnothing\), there are no constraints and
\eqref{eq:general-iso-kkt-proof} reads \(g_u=0\) for all \(u\).
Define
\[
\gamma_*=\max_{(u,v)\in A}\gamma_{uv},
\]
with \(\gamma_*=0\) if \(A=\varnothing\).

Now fix \(\theta\in\R^V\) and write \(\delta=\theta-\theta^*\). The
definition of subgradient gives
\(\ell_u(\theta_u)-\ell_u(\theta^*_u)\geq g_u\delta_u\). Summing this inequality over
\(u\in V\), and using the first line of \eqref{eq:general-iso-kkt-proof}
give
\[
\begin{aligned}
f(\theta)-f^*
&\geq\sum_{u\in V}g_u\delta_u\\
&=\sum_{u\in V}\delta_u\Bigl\{
\sum_{v:(v,u)\in A}\gamma_{vu}
-\sum_{v:(u,v)\in A}\gamma_{uv}\Bigr\}\\
&=\sum_{(u,v)\in A}\gamma_{uv}(\delta_v-\delta_u).
\end{aligned}
\]
The last equality holds because each edge \((u,v)\in A\) appears
exactly twice in the double sum: once with coefficient
\(+\gamma_{uv}\delta_v\), at its head \(v\), and once with
coefficient \(-\gamma_{uv}\delta_u\), at its tail \(u\). By the
second line of \eqref{eq:general-iso-kkt-proof},
\[
\begin{aligned}
\gamma_{uv}(\delta_v-\delta_u)
&=\gamma_{uv}(\theta_v-\theta_u)
-\gamma_{uv}(\theta^*_v-\theta^*_u)\\
&=\gamma_{uv}(\theta_v-\theta_u).
\end{aligned}
\]
Finally, \(x\geq-x_-\) for every \(x\in\R\), and
\(0\leq\gamma_{uv}\leq\gamma_*\), so
\[
\gamma_{uv}(\theta_v-\theta_u)
\geq-\gamma_{uv}(\theta_v-\theta_u)_-
\geq-\gamma_*(\theta_v-\theta_u)_-.
\]
Summing over \((u,v)\in A\) yields
\[
f(\theta)-f^*\geq-\gamma_*P(\theta),
\]
and adding \(\lambda P(\theta)\) to both sides gives
\eqref{eq:exact-penalty-gap}.

\emph{Part (c).} Let \(\lambda>\gamma_*\). By
\eqref{eq:exact-penalty-gap} and \(P\geq0\), we have
\(H_\lambda(\theta)\geq f^*\) for every \(\theta\in\R^V\). If
\(\theta\in\cS_{\mathrm{iso}}\), then \(P(\theta)=0\) by
\eqref{eq:P-zero-iff-isotonic} and \(f(\theta)=f^*\), so
\(H_\lambda(\theta)=f^*\). Thus
\(\min_{\R^V}H_\lambda=f^*\), and every element of
\(\cS_{\mathrm{iso}}\) minimizes \(H_\lambda\). Conversely, let
\(\theta\) minimize \(H_\lambda\), so that \(H_\lambda(\theta)=f^*\).
Then \eqref{eq:exact-penalty-gap} gives
\((\lambda-\gamma_*)P(\theta)\leq0\), and hence \(P(\theta)=0\).
By \eqref{eq:P-zero-iff-isotonic}, \(\theta\in\mathcal K_D\), and
\(f(\theta)=H_\lambda(\theta)=f^*\). Therefore,
\(\theta\in\cS_{\mathrm{iso}}\).
\end{proof}

\begin{proof}[Proof of Theorem~\ref{thm:isotonic}]
We prove the four parts in order. The proof of part~\textup{(ii)}
is organized as the proof of Lemma~\ref{lem:neariso-equals-iso}.

\medskip
\noindent\textbf{Part (i): Coordinate intervals.}

By Lemma~\ref{lem:general-neariso-equals-iso}(a),
\(\cS_{\mathrm{iso}}\) is nonempty, compact, and convex. For each
\(i\in V\), the set \(\{\theta_i:\theta\in\cS_{\mathrm{iso}}\}\) is
the image of \(\cS_{\mathrm{iso}}\) under the continuous linear map
\(\theta\mapsto\theta_i\). It is therefore a nonempty compact convex
subset of \(\R\), that is, a closed bounded interval. Its endpoints
are
\[
\lo_i^{\,\mathrm{iso}}
=\min_{\theta\in\cS_{\mathrm{iso}}}\theta_i,
\qquad
\hi_i^{\,\mathrm{iso}}
=\max_{\theta\in\cS_{\mathrm{iso}}}\theta_i,
\]
both of which are finite. This proves part~\textup{(i)}.

\medskip
\noindent\textbf{Part (ii): The nested-set formula.}

Fix \(i\in V\) and put
\[
t=\hi_i^{\,\mathrm{iso}}\in\R.
\]
The one-sided derivatives of the finite convex functions \(\ell_u\)
are finite at \(t\), so
\[
M:=\max_{u\in V}
\max\bigl\{|\ell'_{u,-}(t)|,\,|\ell'_{u,+}(t)|\bigr\}<\infty.
\]
With \(\gamma_*\) as in
Lemma~\ref{lem:general-neariso-equals-iso}, fix
\begin{equation}\label{eq:iso-general-threshold}
\lambda>\lambda_0:=\max\{\gamma_*,\,|V|M\}.
\end{equation}
This threshold plays the role of
\(\lambda_0=n(y_{\max}-y_{\min})\) in the proof of
Lemma~\ref{lem:neariso-equals-iso}. It may depend on \(i\) through
\(t\). This is harmless: \(\lambda\) is an auxiliary parameter,
which does not appear in \eqref{eq:isotonic-lower} or
\eqref{eq:isotonic-upper}.

By Lemma~\ref{lem:general-neariso-as-tvd}, \(H_\lambda\) is the
generalized graph TVD objective on \(\widetilde G\) with finite
convex losses \(\ell^{(\lambda)}_u\) and penalty parameter
\(\lambda/2\). The tilted losses \(\ell^{(\lambda)}_u\) need not be coercive.
However, \(P\geq0\) by \eqref{eq:P-zero-iff-isotonic}, so
\(H_\lambda\geq f\), and \(f\) is coercive on \(\R^V\) by
\eqref{eq:coercive-on-cone}. Hence \(H_\lambda\) is coercive, and
Theorem~\ref{thm:convex} applies to it. By
Lemma~\ref{lem:general-neariso-equals-iso}(c), the solution set of
this TVD problem is \(\cS_{\mathrm{iso}}\), so its upper endpoint at
coordinate \(i\) is \(t\). Therefore, it follows from
\eqref{eq:convex-upper} that
\begin{equation}\label{eq:iso-general-raw}
t=\min_{\substack{J\subseteq V\\i\in J}}
\max_{\substack{I\subseteq J\\i\in I}}
r_+\bigl(\Phi^{(\lambda)}_{I,J}\bigr),
\end{equation}
where
\[
\Phi^{(\lambda)}_{I,J}(x)
=\sum_{u\in I}\ell^{(\lambda)}_u(x)
+\frac{\lambda}{2}\,C_{I,J}\,x,
\]
and \(C_{I,J}\) is computed in \(\widetilde G\).

As in Section~\ref{sec:proof-isotonic-as-tvd}, for disjoint
\(S,T\subseteq V\) let \(A(S,T)\) denote the set of directed edges
of \(D\) from \(S\) to \(T\). By the definition
\eqref{eq:tilted-losses} of \(\ell^{(\lambda)}_u\) and the boundary
identity \eqref{eq:iso-proof-boundary},
\[
\begin{aligned}
\Phi^{(\lambda)}_{I,J}(x)
&=\sum_{u\in I}\ell_u(x)
+\frac{\lambda}{2}\Bigl\{C_{I,J}+\sum_{u\in I}d_u\Bigr\}x\\
&=\sum_{u\in I}\ell_u(x)\\
&\quad+\lambda\bigl\{|A(I,J\setminus I)|-|A(J^c,I)|\bigr\}x.
\end{aligned}
\]
Accordingly, for \(i\in I\subseteq J\subseteq V\), define
\begin{equation}\label{eq:iso-general-Q}
\begin{aligned}
\Psi_{I,J}(x)&=\sum_{u\in I}\ell_u(x)\\
&\quad+\lambda\bigl\{|A(I,J\setminus I)|-|A(J^c,I)|\bigr\}x,\\
Q_{I,J}&=r_+(\Psi_{I,J})\in[-\infty,+\infty],\\
\mu_J&=\max_{\substack{I\subseteq J\\i\in I}}Q_{I,J}.
\end{aligned}
\end{equation}
With this notation, \eqref{eq:iso-general-raw} becomes
\begin{equation}\label{eq:iso-general-nested}
t=\min_{\substack{J\subseteq V\\i\in J}}\mu_J
=\min_{\substack{J\subseteq V\\i\in J}}
\max_{\substack{I\subseteq J\\i\in I}}Q_{I,J}.
\end{equation}
This is the counterpart of \eqref{eq:iso-proof-nested}: under
squared loss,
\(Q_{I,J}=\bar y_I
-\lambda\{|A(I,J\setminus I)|-|A(J^c,I)|\}/|I|\).
The minima and maxima in \eqref{eq:iso-general-nested} are taken
over finite collections of extended real numbers, so they are
attained. Individual values \(Q_{I,J}\) may be infinite, but the
left-hand side \(t\) is finite.

In the squared-loss proof, \eqref{eq:iso-proof-nested} was
simplified by comparing its terms with \(y_{\max}\) and
\(y_{\min}\). Here we compare the terms \(Q_{I,J}\) with the value
\(t\) itself, through the sign of a one-sided derivative of
\(\Psi_{I,J}\) at \(t\) and
Lemma~\ref{lem:endpoint-comparisons}. Since one-sided derivatives
are additive, for every \(i\in I\subseteq J\),
\begin{equation}\label{eq:iso-general-derivatives}
\begin{aligned}
\Psi'_{I,J,\pm}(t)
&=\sum_{u\in I}\ell'_{u,\pm}(t)\\
&\quad+\lambda\bigl\{|A(I,J\setminus I)|-|A(J^c,I)|\bigr\},
\end{aligned}
\end{equation}
and, by the definition of \(M\),
\begin{equation}\label{eq:iso-general-M-bound}
\Bigl|\sum_{u\in I}\ell'_{u,\pm}(t)\Bigr|
\leq|I|M\leq|V|M<\lambda.
\end{equation}

\medskip
\noindent\textbf{Step 1: Restricting the outer minimum to lower sets.}

Let \(J\ni i\) be a set that is not lower. Then some edge of \(D\)
enters \(J\) from outside, so
\[
k:=|A(J^c,J)|\geq1.
\]
The choice \(I=J\) is admissible in the inner maximum, and
\(A(J,J\setminus J)=\varnothing\). By
\eqref{eq:iso-general-derivatives} and
\eqref{eq:iso-general-M-bound},
\[
\Psi'_{J,J,+}(t)
=\sum_{u\in J}\ell'_{u,+}(t)-\lambda k
\leq|V|M-\lambda<0.
\]
The fourth assertion of Lemma~\ref{lem:endpoint-comparisons} gives
\(Q_{J,J}>t\). Consequently,
\begin{equation}\label{eq:iso-general-nonlower}
\mu_J\geq Q_{J,J}>t
\quad\text{for every non-lower }J\ni i.
\end{equation}
By \eqref{eq:iso-general-nested}, \(t\) is the smallest of the
values \(\mu_J\), and this minimum is attained. By
\eqref{eq:iso-general-nonlower}, it is not attained at any
non-lower set. Since \(V\in\cL_i\), the collection \(\cL_i\) is
nonempty, and therefore
\[
t=\min_{L\in\cL_i}\mu_L.
\]
Now let \(L\in\cL_i\) and \(i\in I\subseteq L\). No edge of \(D\)
enters \(L\) from \(L^c\), so \(A(L^c,I)=\varnothing\) and
\[
\Psi_{I,L}(x)
=\sum_{u\in I}\ell_u(x)+\lambda|A(I,L\setminus I)|\,x.
\]
We have thus shown that
\begin{equation}\label{eq:iso-general-lower}
\begin{aligned}
t&=\min_{L\in\cL_i}
\max_{\substack{I\subseteq L\\i\in I}}Q_{I,L},\\
Q_{I,L}
&=r_+\Bigl(x\mapsto\sum_{u\in I}\ell_u(x)
+\lambda|A(I,L\setminus I)|\,x\Bigr),
\end{aligned}
\end{equation}
which is the counterpart of \eqref{eq:iso-proof-lower}.

\medskip
\noindent\textbf{Step 2: Expressing the inner maximum using upper sets.}

Fix \(L\in\cL_i\). By \eqref{eq:iso-general-lower},
\begin{equation}\label{eq:iso-general-mu-lower}
\mu_L\geq t.
\end{equation}
Let \(i\in I\subseteq L\) be a set with
\[
k:=|A(I,L\setminus I)|\geq1.
\]
By \eqref{eq:iso-general-derivatives} with
\(A(L^c,I)=\varnothing\), and by \eqref{eq:iso-general-M-bound},
\[
\Psi'_{I,L,-}(t)
=\sum_{u\in I}\ell'_{u,-}(t)+\lambda k
\geq-|V|M+\lambda>0.
\]
The fifth assertion of Lemma~\ref{lem:endpoint-comparisons} and
\eqref{eq:iso-general-mu-lower} give
\[
Q_{I,L}<t\leq\mu_L.
\]
Hence no set \(I\) with \(A(I,L\setminus I)\neq\varnothing\)
attains the inner maximum \(\mu_L\). The set \(I=L\) satisfies
\(A(L,L\setminus L)=\varnothing\), so the remaining collection is
nonempty, and
\[
\mu_L=\max\bigl\{Q_{I,L}:
i\in I\subseteq L,\ A(I,L\setminus I)=\varnothing\bigr\}.
\]
For every such \(I\), the linear term of \(\Psi_{I,L}\) vanishes,
so that, by \eqref{eq:isotonic-scalars},
\[
Q_{I,L}
=r_+\Bigl(x\mapsto\sum_{u\in I}\ell_u(x)\Bigr)
=q^+_I.
\]

It remains to identify the surviving sets. This is a statement
about the DAG \(D\) alone, and it was established in Step~2 of the
proof of Lemma~\ref{lem:neariso-equals-iso}:
\[
\begin{aligned}
&\bigl\{I:i\in I\subseteq L,\
A(I,L\setminus I)=\varnothing\bigr\}\\
&\qquad=\{L\cap U:U\in\cU_i\}.
\end{aligned}
\]
Therefore
\[
\mu_L=\max_{U\in\cU_i}q^+_{L\cap U}
\qquad\text{for every }L\in\cL_i,
\]
and substituting this into \(t=\min_{L\in\cL_i}\mu_L\) gives
\begin{equation}\label{eq:iso-general-minmax}
\hi_i^{\,\mathrm{iso}}
=\min_{L\in\cL_i}\max_{U\in\cU_i}q^+_{L\cap U}.
\end{equation}
This is \eqref{eq:isotonic-upper}, and the counterpart of
\eqref{eq:iso-proof-minmax}. The right-hand side is finite, since
the left-hand side is.

\medskip
\noindent\textbf{Step 3: The maxmin identity.}

The preceding argument used only the standing assumptions. Hence
\eqref{eq:iso-general-minmax} holds for every vertex of every
finite DAG and every family of finite, convex, coercive losses. We apply it to a
reflected problem, as in Step~3 of the proof of
Lemma~\ref{lem:neariso-equals-iso}.

Let \(D^{\mathrm{rev}}=(V,A^{\mathrm{rev}})\), with
\(A^{\mathrm{rev}}=\{(v,u):(u,v)\in A\}\), be obtained by
reversing every edge of \(D\), and let
\[
\widetilde\ell_u(x)=\ell_u(-x),
\qquad
\widetilde f(\theta)=\sum_{u\in V}\widetilde\ell_u(\theta_u)
=f(-\theta).
\]
The reflected losses remain finite, convex, and coercive. Since
\(\mathcal K_{D^{\mathrm{rev}}}=-\mathcal K_D\), the reflected
solution set is \(-\cS_{\mathrm{iso}}\), with upper endpoint
\(-\lo_i^{\,\mathrm{iso}}\) at coordinate \(i\). Reversing the
edges interchanges lower and upper sets. Moreover, for nonempty
\(B\subseteq V\), the sixth assertion of
Lemma~\ref{lem:endpoint-comparisons} gives
\[
r_+\Bigl(x\mapsto\sum_{u\in B}\widetilde\ell_u(x)\Bigr)
=-r_-\Bigl(\sum_{u\in B}\ell_u\Bigr)
=-q^-_B.
\]
Applying \eqref{eq:iso-general-minmax} to the reflected problem
therefore gives
\[
-\lo_i^{\,\mathrm{iso}}
=\min_{U\in\cU_i}\max_{L\in\cL_i}\bigl(-q^-_{L\cap U}\bigr).
\]
Negating both sides exchanges minimum and maximum:
\[
\lo_i^{\,\mathrm{iso}}
=\max_{U\in\cU_i}\min_{L\in\cL_i}q^-_{L\cap U}.
\]
This is \eqref{eq:isotonic-lower}, which completes the proof of
part~\textup{(ii)}.

\medskip
\noindent\textbf{Part (iii): Lattice structure and simultaneous
attainment.}

Let \(\theta,\eta\in\cS_{\mathrm{iso}}\). We first show that
\(\theta\wedge\eta\) and \(\theta\vee\eta\) belong to
\(\mathcal K_D\). For every \((u,v)\in A\), we have
\(\theta_u\leq\theta_v\) and \(\eta_u\leq\eta_v\), and therefore
\[
\begin{aligned}
\min\{\theta_u,\eta_u\}&\leq\min\{\theta_v,\eta_v\},\\
\max\{\theta_u,\eta_u\}&\leq\max\{\theta_v,\eta_v\}.
\end{aligned}
\]
Next, for every \(u\in V\), the pair
\((\min\{\theta_u,\eta_u\},\max\{\theta_u,\eta_u\})\) is a
rearrangement of \((\theta_u,\eta_u)\). Summing
\(\ell_u\) over both entries and then over \(u\in V\) gives
\[
f(\theta\wedge\eta)+f(\theta\vee\eta)
=f(\theta)+f(\eta)=2f^*.
\]
Since both \(\theta\wedge\eta\) and \(\theta\vee\eta\) lie in
\(\mathcal K_D\), each term on the left is at least \(f^*\).
Hence both terms equal \(f^*\), that is,
\(\theta\wedge\eta,\theta\vee\eta\in\cS_{\mathrm{iso}}\). This
proves that \(\cS_{\mathrm{iso}}\) is a lattice.

By part~\textup{(i)}, for each \(j\in V\) we can choose
\(\theta^{j,-},\theta^{j,+}\in\cS_{\mathrm{iso}}\) with
\[
\theta^{j,-}_j=\lo_j^{\,\mathrm{iso}},
\qquad
\theta^{j,+}_j=\hi_j^{\,\mathrm{iso}}.
\]
Since \(V\) is finite, applying the lattice property finitely many
times gives
\[
\theta^-:=\bigwedge_{j\in V}\theta^{j,-}\in\cS_{\mathrm{iso}},
\qquad
\theta^+:=\bigvee_{j\in V}\theta^{j,+}\in\cS_{\mathrm{iso}}.
\]
Fix \(u\in V\). On the one hand,
\(\theta^-_u\leq\theta^{u,-}_u=\lo_u^{\,\mathrm{iso}}\). On the
other hand, \(\theta^-\in\cS_{\mathrm{iso}}\), so
\(\theta^-_u\geq\lo_u^{\,\mathrm{iso}}\) by the definition of
\(\lo_u^{\,\mathrm{iso}}\). Thus
\(\theta^-=\lo^{\,\mathrm{iso}}\), and in the same way
\(\theta^+=\hi^{\,\mathrm{iso}}\). Consequently,
\(\lo^{\,\mathrm{iso}}\) and \(\hi^{\,\mathrm{iso}}\) belong to
\(\cS_{\mathrm{iso}}\), and by the definition of the coordinate
endpoints every \(\theta\in\cS_{\mathrm{iso}}\) satisfies
\(\lo_u^{\,\mathrm{iso}}\leq\theta_u\leq\hi_u^{\,\mathrm{iso}}\)
for all \(u\in V\). They are therefore the coordinatewise smallest
and largest solutions, which proves part~\textup{(iii)}.

\medskip
\noindent\textbf{Part (iv): Uniqueness under strict convexity.}

Suppose that \(f\) is strictly convex on \(\mathcal K_D\), and let
\(\theta\neq\eta\) be two elements of \(\cS_{\mathrm{iso}}\). Their
midpoint \((\theta+\eta)/2\) belongs to the convex set
\(\mathcal K_D\), and strict convexity gives
\[
f\Bigl(\frac{\theta+\eta}{2}\Bigr)
<\frac{f(\theta)+f(\eta)}{2}=f^*,
\]
contradicting the definition of \(f^*\). Hence
\(\cS_{\mathrm{iso}}\), which is nonempty by part~\textup{(i)},
consists of a single vector, and
\(\lo^{\,\mathrm{iso}}=\hi^{\,\mathrm{iso}}\) by
part~\textup{(iii)}. Finally, if every \(\ell_u\) is strictly
convex, then \(f\) is strictly convex on \(\R^V\): for
\(\theta\neq\eta\), the strict inequality holds in at least one
coordinate and the weak inequality in all others. In particular,
\(f\) is strictly convex on \(\mathcal K_D\). This proves
part~\textup{(iv)}.
\end{proof}

\section{Proofs of Auxiliary Facts}\label{sec:proof_aux}

\subsection{Proof of Remark~\ref{rem:convex_fn_min_endpoint}}\label{subsec:proof-convex-endpoint}
Suppose the set of minimizers of $h$ is non-empty. Due to convexity, the set of minimizers of $h$ has to be a convex set in $\mathbb R$. To see this, note that for any minimizers $x_1\leq x_2$ and any $\alpha\in[0,1]$, convexity gives
\[h\left(\alpha x_1+(1-\alpha)x_2\right)\leq \alpha h(x_1)+(1-\alpha)h(x_2)=\min_{x\in \mathbb R} h(x).\]
Now, any convex set in $\mathbb R$ is an interval, possibly unbounded or a singleton. Moreover, the set of minimizers is closed due to continuity of $h$. This proves first part of the remark.\\
\\
Now, suppose the set of minimizers of $h$ is bounded. The minimizer set is also closed, so it is compact. Let us denote this interval as $[m,M]$. Now, it holds that $x\in \mathbb R$ minimizes the convex map $h$, if and only if
\[h'_-(x)\leq 0 \leq h'_+(x).\]
This follows easily from the subdifferential characterization in \eqref{eqn:charac_subgradient} and the associated optimality conditions. Therefore $h'_-(m)\leq 0\leq h'_+(m)$, and for any $\epsilon>0$, due to monotonicity of one-sided derivative of convex map $h$, we have 
\[h'_-(m-\epsilon)\leq h'_-(m)\leq 0.\]
If $h'_+(m-\epsilon)\geq 0$, then it would imply that $m-\epsilon$ also minimizes $h$, which contradicts the definition of $m$. Therefore in this case,
\[m=\min\{t:h'_+(t)\geq 0\}=r_-(h).\]
It can be similarly proved that $M=\max\{t:h'_-(t)\leq 0\}=r_+(h).$ 
\hfill$\square$

\subsection{Proof of Remark~\ref{rem:logistic}}
\label{subsec:proof-logistic}

Suppose \(F_\lambda(\theta)\le M\), where \(M\ge0\). Since all
loss terms and the TV penalty are nonnegative, each loss is at
most \(M\), and \(\TV_G(\theta)\le M/\lambda\). For any vertices
\(u,v\), connectedness provides a simple path from \(u\) to \(v\).
Applying the triangle inequality along this path gives
\[
|\theta_u-\theta_v|
\le\TV_G(\theta)\le M/\lambda.
\]
Thus the TV penalty bounds the difference between any two
coordinates.

Choose vertices \(a,b\) with \(y_a=0\) and \(y_b=1\). The loss at
\(a\) satisfies
\[
M\ge\ell_a(\theta_a)
=\log(1+e^{\theta_a})\ge\theta_a,
\]
so \(\theta_a\le M\). Consequently, for every \(u\in V\),
\[
\theta_u\le\theta_a+|\theta_u-\theta_a|
\le M+M/\lambda.
\]
Similarly, the loss at \(b\) satisfies
\[
M\ge\ell_b(\theta_b)
=\log(1+e^{-\theta_b})\ge-\theta_b,
\]
so \(\theta_b\ge-M\), and hence
\[
\theta_u\ge\theta_b-|\theta_u-\theta_b|
\ge-M-M/\lambda.
\]
Therefore every sublevel set of \(F_\lambda\) is bounded, proving
coercivity. The losses are finite and convex, so
Theorem~\ref{thm:convex} applies.

Finally, each logistic loss is strictly convex, since
\[
\ell_u''(t)=\frac{e^t}{(1+e^t)^2}>0
\qquad\text{for every }t\in\R.
\]
Their separable sum is therefore strictly convex on \(\R^V\).
Adding the convex TV penalty preserves strict convexity, so
\(F_\lambda\) has a unique minimizer.
\hfill\(\square\)

\subsection{Proof of Lemma~\ref{lem:convex-function-facts}}
Recall from \eqref{eq:scalar-endpoints} that
\[
\begin{aligned}
r_+(h)&=\sup\{x:h'_-(x)\leq0\},\\
r_-(h)&=\inf\{x:h'_+(x)\geq0\}.
\end{aligned}
\]
By \citet[Theorem~24.1]{rockafellar1997convex}, \(h'_-\) and
\(h'_+\) are finite and nondecreasing on \(\R\), and
\begin{equation}\label{eq:one-sided-limits}
\lim_{z\downarrow t}h'_-(z)=h'_+(t),
\qquad
\lim_{z\uparrow t}h'_-(z)=h'_-(t).
\end{equation}

\emph{1.} If \(h'_-(t)\leq0\), then \(t\) belongs to
the set \(\{x:h'_-(x)\leq0\}\), whose supremum is \(r_+(h)\).

\emph{2.} If \(r_+(h)>t\), the definition of the
supremum gives a point \(s>t\) with \(h'_-(s)\leq0\). Fix any
\(0<\varepsilon<s-t\). 
We will use the secant
inequality for real-valued convex maps: for \(a<b<c\),
\[
\frac{h(b)-h(a)}{b-a}
\leq\frac{h(c)-h(a)}{c-a}
\leq\frac{h(c)-h(b)}{c-b}.
\]
It follows from the secant-inequality that for any $0<\delta<\epsilon$,
\[\frac{h(t+\varepsilon)-h(t)}{\varepsilon}\leq \frac{h(t+\varepsilon)-h(t+\delta)}{\varepsilon-\delta}.\]
Taking $\delta\uparrow \varepsilon$, the RHS converges to $h'_-(t+\varepsilon)$. Therefore, using the monotonicity of \(h'_-\), we get
\[
\frac{h(t+\varepsilon)-h(t)}{\varepsilon}
\leq h'_-(t+\varepsilon)\leq h'_-(s)\leq0.
\]

\emph{3.} \ Fix \(\varepsilon>0\). We need to show that
\(h(x+\varepsilon)-h(x)\leq h(y+\varepsilon)-h(y)\) whenever
\(x\leq y\), the case \(x=y\) being trivial. 
First, consider the case \(x<y<x+\varepsilon\). Applying the
secant inequality to \(x<y<x+\varepsilon\), and then to
\(y<x+\varepsilon<y+\varepsilon\), gives
\[
\frac{h(x+\varepsilon)-h(x)}{\varepsilon}
\leq\frac{h(x+\varepsilon)-h(y)}{x+\varepsilon-y}
\leq\frac{h(y+\varepsilon)-h(y)}{\varepsilon}.
\]

Next, consider the case \(x+\varepsilon<y\). Applying the secant
inequality to \(x<x+\varepsilon<y\), and then to
\(x+\varepsilon<y<y+\varepsilon\), gives
\[
\frac{h(x+\varepsilon)-h(x)}{\varepsilon}
\leq\frac{h(y)-h(x+\varepsilon)}{y-x-\varepsilon}
\leq\frac{h(y+\varepsilon)-h(y)}{\varepsilon}.
\]

Finally, if \(y=x+\varepsilon\), the secant inequality applied to
\(x<x+\varepsilon<x+2\varepsilon\) gives the claim directly.

In every case the two extreme terms yield
\[
h(x+\varepsilon)-h(x)\leq h(y+\varepsilon)-h(y),
\]
as desired.

\emph{4.} If \(h'_+(t)<0\), the first limit in
\eqref{eq:one-sided-limits} gives a point \(s>t\) with \(h'_-(s)<0\).
Then \(s\in\{x:h'_-(x)\leq0\}\), so \(r_+(h)\geq s>t\).

\emph{5.} If \(h'_-(t)>0\), the second limit in
\eqref{eq:one-sided-limits} gives a point \(s<t\) with \(h'_-(s)>0\).
Since \(h'_-\) is nondecreasing, \(h'_-(x)>0\) for every \(x\geq s\).
Thus \(\{x:h'_-(x)\leq0\}\subseteq(-\infty,s)\), and
\(r_+(h)\leq s<t\); if the set is empty, \(r_+(h)=-\infty<t\).

\emph{6.} For every \(x\in\R\),
\[
\widetilde h'_-(x)
=\lim_{\varepsilon\downarrow0}
\frac{h(-x)-h(-x+\varepsilon)}{\varepsilon}
=-h'_+(-x).
\]
Therefore, substituting \(s=-x\),
\[
\begin{aligned}
r_+(\widetilde h)
&=\sup\{x:h'_+(-x)\geq0\}\\
&=-\inf\{s:h'_+(s)\geq0\}
=-r_-(h),
\end{aligned}
\]
where the conventions \(\sup\varnothing=-\infty\) and
\(\inf\varnothing=+\infty\) are consistent with each other.
\hfill$\square$

\medskip
\section*{Acknowledgement}
Sabyasachi Chatterjee thanks Ashwin Pananjady for a conversation which initiated this research.

\bibliographystyle{plainnat}
\bibliography{TVD_References}
\end{document}